\input amstex
\input amsppt.sty

\input epsf
\input texdraw

\magnification1200
\hsize14cm
\vsize19cm

\TagsOnRight

%figure numbers

\def\FB{5}
\def\FC{2}
\def\FD{7}
\def\FE{3}
\def\FF{4}
\def\FFa{1}

%equation numbers

\def\Bmm{2.1}
\def\Bm{2.2}
\def\TABzero{2.3}
\def\TABone{2.4}
\def\Bp{2.5}
\def\Bpp{2.6}

\def\SPP{3.1}
\def\AAAAAA{4.1}
\def\PPF{4.2}
\def\PF{4.3}
\def\AAAAA{4.4}
\def\AAAAAAA{4.5}
\def\AAAA{4.6}
\def\AAAa{4.7}
\def\AAAb{4.8}
\def\AAa{4.9}
\def\AA{4.10}
\def\AAb{4.11}
\def\AB{4.12}
\def\ABb{4.13}
\def\ABa{3.2}
\def\SS{5.2}
\def\ABB{5.1}
\def\BD{5.3}
\def\AC{5.4}
\def\AD{5.5}
\def\AE{5.6}
\def\AFa{5.7}
\def\AFb{5.8}
\def\AF{5.9}
\def\AG{5.10}
\def\AH{5.11}
\def\BAa{5.12}
\def\BA{5.13}
\def\BB{5.14}
\def\BCa{5.15}
\def\BCb{5.16}
\def\BCc{5.17}
\def\BCd{5.18}
\def\BCe{5.19}
\def\BCf{5.20}
\def\BCg{5.21}
\def\HYP{5.22}
\def\BCh{5.23}
\def\BCi{5.24}
\def\BF{5.25}
\def\BGa{5.26}
\def\BG{5.27}
\def\BH{5.28}
\def\BI{5.29}
\def\CAa{6.1}
\def\CAb{6.2}
\def\CAc{6.3}
\def\CAe{6.4}
\def\CAd{6.5}
\def\CAg{6.6}
\def\CAf{6.7}
\def\CA{7.1}
\def\CB{7.7}
\def\CC{7.8}
\def\CD{7.9}
\def\CE{7.10}
\def\X{7.11}
\def\CF{7.12}
\def\CG{7.13}
\def\CH{7.14}
\def\CI{7.15}
\def\CJ{7.16}
\def\CK{7.17}
\def\CL{7.18}
\def\CM{7.19}

%theorem numbers
\def\TAA{1}
\def\TAB{2}
\def\TAC{3}
\def\TAD{5}
\def\TAE{6}
\def\TAF{7}
\def\TAG{8}
\def\TB{4}
\def\TCa{9}
\def\TCb{10}
\def\TCc{11}
\def\TC{12}
\def\TD{13}
\def\TE{14}
\def\TF{15}

%reference numbers
\def\AndrAK{1}
\def\AnARAA{2}
\def\BailAA{3}
\def\CiucAK{4}

\def\sc{5}
\def\CiucAM{6}
\def\ec{6}
\def\CiucAO{7}
\def\ov{7}
\def\ef{8}
\def\CLP{9}
\def\FrReAA{10}
\def\Feytwo{11}
\def\FiscAH{12}
\def\FS{13}
\def\GeViAB{14}
\def\GordAC{15}
\def\GrKPAA{16}
\def\IsWaAA{17}
\def\Kone{18}
\def\Kmonbd{19}
\def\KOS{20}
\def\KratAP{21}
\def\KratBI{22}
\def\KratBN{23}
\def\LindAA{24}
\def\MacdAC{25}
\def\MacMAA{26}
\def\MeWaAA{27}
\def\OlveAA{28}
\def\ProcAD{29}
\def\SlatAC{30}
\def\StemAE{31}
\def\VM{32}

\def\po#1#2{(#1)_#2}
\def\({\left(}
\def\){\right)}
\def\fl#1{\lfloor#1\rfloor}

\def\Ga{\Gamma}
\def\al{\alpha}
\def\be{\beta}
\def\ep{\varepsilon}
\def\om{\omega}
\def\sgn{\operatorname{sgn}}
\def\Pf{\operatorname{Pf}}
\def\M{\operatorname{M}}
\def\de{\operatorname{d}}
\def\Z{\Bbb Z}

%Christian's texdraw definitions

\def\rdreieck{\bsegment
  \rlvec(0.866025403784439 -.5) \rlvec(-0.866025403784439 -.5)  \rlvec(0 1)
  \savepos(0 -1)(*ex *ey)
        \esegment
  \move(*ex *ey)
        }
\def\rhombus{\bsegment
  \rlvec(0.866025403784439 .5) \rlvec(0.866025403784439 -.5) 
  \rlvec(-0.866025403784439 -.5)  \rlvec(0 1)        
  \rmove(0 -1)  \rlvec(-0.866025403784439 .5) 
  \savepos(0.866025403784439 -.5)(*ex *ey)
        \esegment
  \move(*ex *ey)
        }
\def\RhombusA{\bsegment
  \rlvec(0.866025403784439 .5) \rlvec(0.866025403784439 -.5) 
  \rlvec(-0.866025403784439 -.5) \rlvec(-0.866025403784439 .5) 
  \savepos(0.866025403784439 -.5)(*ex *ey)
        \esegment
  \move(*ex *ey)
        }
\def\RhombusB{\bsegment
  \rlvec(0.866025403784439 .5) \rlvec(0 -1)
  \rlvec(-0.866025403784439 -.5) \rlvec(0 1) 
  \savepos(0 -1)(*ex *ey)
        \esegment
  \move(*ex *ey)
        }
\def\RhombusC{\bsegment
  \rlvec(0.866025403784439 -.5) \rlvec(0 -1)
  \rlvec(-0.866025403784439 .5) \rlvec(0 1) 
  \savepos(0.866025403784439 -.5)(*ex *ey)
        \esegment
  \move(*ex *ey)
        }

\def\huSchritt{\bsegment
  \rlvec(0.866025403784439 -.5) 
  \savepos(0.866025403784439 -.5)(*ex *ey)
        \esegment
  \move(*ex *ey)
        }
\def\hoSchritt{\bsegment
  \rlvec(0.866025403784439 .5) 
  \savepos(0.866025403784439 .5)(*ex *ey)
        \esegment
  \move(*ex *ey)
        }

\def\vdSchritt{\bsegment
  \lpatt(.05 .13)
  \rlvec(0 -1) 
  \savepos(0 -1)(*ex *ey)
        \esegment
  \move(*ex *ey)
        }

\def\odSchritt{\bsegment
  \lpatt(.05 .13)
  \rlvec(-0.866025403784439 -.5) 
  \savepos(-0.866025403784439 -.5)(*ex *ey)
        \esegment
  \move(*ex *ey)
        }

\def\ringerl(#1 #2){\move(#1 #2)\fcir f:0 r:.15}
\def\Ringerl(#1 #2){\move(#1 #2)\lcir r:.19}

\catcode`\@=11
\font\tenln    = line10
\font\tenlnw   = linew10

\newskip\Einheit \Einheit=0.5cm
\newcount\xcoord \newcount\ycoord
\newdimen\xdim \newdimen\ydim \newdimen\PfadD@cke \newdimen\Pfadd@cke

%%%%%%%%%%%%%%%%%%%%%%%%%%%%%%%%%%%%%%%%%%%%%%%%%
%LaTeX counters, dimensions, variables for lines%
%%%%%%%%%%%%%%%%%%%%%%%%%%%%%%%%%%%%%%%%%%%%%%%%%
\newcount\@tempcnta
\newcount\@tempcntb

\newdimen\@tempdima
\newdimen\@tempdimb

\newdimen\@wholewidth
\newdimen\@halfwidth

\newcount\@xarg
\newcount\@yarg
\newcount\@yyarg
\newbox\@linechar
\newbox\@tempboxa
\newdimen\@linelen
\newdimen\@clnwd
\newdimen\@clnht

\newif\if@negarg

\def\@whilenoop#1{}
\def\@whiledim#1\do #2{\ifdim #1\relax#2\@iwhiledim{#1\relax#2}\fi}
\def\@iwhiledim#1{\ifdim #1\let\@nextwhile=\@iwhiledim
        \else\let\@nextwhile=\@whilenoop\fi\@nextwhile{#1}}

\def\@whileswnoop#1\fi{}
\def\@whilesw#1\fi#2{#1#2\@iwhilesw{#1#2}\fi\fi}
\def\@iwhilesw#1\fi{#1\let\@nextwhile=\@iwhilesw
         \else\let\@nextwhile=\@whileswnoop\fi\@nextwhile{#1}\fi}

\def\thinlines{\let\@linefnt\tenln \let\@circlefnt\tencirc
  \@wholewidth\fontdimen8\tenln \@halfwidth .5\@wholewidth}
\def\thicklines{\let\@linefnt\tenlnw \let\@circlefnt\tencircw
  \@wholewidth\fontdimen8\tenlnw \@halfwidth .5\@wholewidth}
\thinlines
%%%%%%%%%%%%%%%%%%%%%%%%%%%%%%%%%%%%%%%%%%%%%%%%%%%%%%%%%%%

\PfadD@cke1pt \Pfadd@cke0.5pt
\def\PfadDicke#1{\PfadD@cke#1 \divide\PfadD@cke by2 \Pfadd@cke\PfadD@cke \multiply\PfadD@cke by2}
\long\def\LOOP#1\REPEAT{\def\BODY{#1}\ITERATE}
\def\ITERATE{\BODY \let\next\ITERATE \else\let\next\relax\fi \next}
\let\REPEAT=\fi
\def\Punkt{\hbox{\raise-2pt\hbox to0pt{\hss$\ssize\bullet$\hss}}}
\def\DuennPunkt(#1,#2){\unskip
  \raise#2 \Einheit\hbox to0pt{\hskip#1 \Einheit
          \raise-2.5pt\hbox to0pt{\hss$\bullet$\hss}\hss}}
\def\NormalPunkt(#1,#2){\unskip
  \raise#2 \Einheit\hbox to0pt{\hskip#1 \Einheit
          \raise-3pt\hbox to0pt{\hss\twelvepoint$\bullet$\hss}\hss}}
\def\DickPunkt(#1,#2){\unskip
  \raise#2 \Einheit\hbox to0pt{\hskip#1 \Einheit
          \raise-4pt\hbox to0pt{\hss\fourteenpoint$\bullet$\hss}\hss}}
\def\Kreis(#1,#2){\unskip
  \raise#2 \Einheit\hbox to0pt{\hskip#1 \Einheit
          \raise-4pt\hbox to0pt{\hss\fourteenpoint$\circ$\hss}\hss}}

%%%%%%%%%%%%%%%%%%%%%
%LaTeX line macros%
%%%%%%%%%%%%%%%%%%%%%
\def\Line@(#1,#2)#3{\@xarg #1\relax \@yarg #2\relax
\@linelen=#3\Einheit
\ifnum\@xarg =0 \@vline
  \else \ifnum\@yarg =0 \@hline \else \@sline\fi
\fi}

\def\@sline{\ifnum\@xarg< 0 \@negargtrue \@xarg -\@xarg \@yyarg -\@yarg
  \else \@negargfalse \@yyarg \@yarg \fi
\ifnum \@yyarg >0 \@tempcnta\@yyarg \else \@tempcnta -\@yyarg \fi
\ifnum\@tempcnta>6 \@badlinearg\@tempcnta0 \fi
\ifnum\@xarg>6 \@badlinearg\@xarg 1 \fi
\setbox\@linechar\hbox{\@linefnt\@getlinechar(\@xarg,\@yyarg)}%
\ifnum \@yarg >0 \let\@upordown\raise \@clnht\z@
   \else\let\@upordown\lower \@clnht \ht\@linechar\fi
\@clnwd=\wd\@linechar
\if@negarg \hskip -\wd\@linechar \def\@tempa{\hskip -2\wd\@linechar}\else
     \let\@tempa\relax \fi
\@whiledim \@clnwd <\@linelen \do
  {\@upordown\@clnht\copy\@linechar
   \@tempa
   \advance\@clnht \ht\@linechar
   \advance\@clnwd \wd\@linechar}%
\advance\@clnht -\ht\@linechar
\advance\@clnwd -\wd\@linechar
\@tempdima\@linelen\advance\@tempdima -\@clnwd
\@tempdimb\@tempdima\advance\@tempdimb -\wd\@linechar
\if@negarg \hskip -\@tempdimb \else \hskip \@tempdimb \fi
\multiply\@tempdima \@m
\@tempcnta \@tempdima \@tempdima \wd\@linechar \divide\@tempcnta \@tempdima
\@tempdima \ht\@linechar \multiply\@tempdima \@tempcnta
\divide\@tempdima \@m
\advance\@clnht \@tempdima
\ifdim \@linelen <\wd\@linechar
   \hskip \wd\@linechar
  \else\@upordown\@clnht\copy\@linechar\fi}

\def\@hline{\ifnum \@xarg <0 \hskip -\@linelen \fi
\vrule height\Pfadd@cke width \@linelen depth\Pfadd@cke
\ifnum \@xarg <0 \hskip -\@linelen \fi}

\def\@getlinechar(#1,#2){\@tempcnta#1\relax\multiply\@tempcnta 8
\advance\@tempcnta -9 \ifnum #2>0 \advance\@tempcnta #2\relax\else
\advance\@tempcnta -#2\relax\advance\@tempcnta 64 \fi
\char\@tempcnta}

\def\Vektor(#1,#2)#3(#4,#5){\unskip\leavevmode
  \xcoord#4\relax \ycoord#5\relax
      \raise\ycoord \Einheit\hbox to0pt{\hskip\xcoord \Einheit
         \Vector@(#1,#2){#3}\hss}}

\def\Vector@(#1,#2)#3{\@xarg #1\relax \@yarg #2\relax
\@tempcnta \ifnum\@xarg<0 -\@xarg\else\@xarg\fi
\ifnum\@tempcnta<5\relax
\@linelen=#3\Einheit
\ifnum\@xarg =0 \@vvector
  \else \ifnum\@yarg =0 \@hvector \else \@svector\fi
\fi
\else\@badlinearg\fi}

\def\@hvector{\@hline\hbox to 0pt{\@linefnt
\ifnum \@xarg <0 \@getlarrow(1,0)\hss\else
    \hss\@getrarrow(1,0)\fi}}

\def\@vvector{\ifnum \@yarg <0 \@downvector \else \@upvector \fi}

\def\@svector{\@sline
\@tempcnta\@yarg \ifnum\@tempcnta <0 \@tempcnta=-\@tempcnta\fi
\ifnum\@tempcnta <5
  \hskip -\wd\@linechar
  \@upordown\@clnht \hbox{\@linefnt  \if@negarg
  \@getlarrow(\@xarg,\@yyarg) \else \@getrarrow(\@xarg,\@yyarg) \fi}%
\else\@badlinearg\fi}

\def\@upline{\hbox to \z@{\hskip -.5\Pfadd@cke \vrule width \Pfadd@cke
   height \@linelen depth \z@\hss}}

\def\@downline{\hbox to \z@{\hskip -.5\Pfadd@cke \vrule width \Pfadd@cke
   height \z@ depth \@linelen \hss}}

\def\@upvector{\@upline\setbox\@tempboxa\hbox{\@linefnt\char'66}\raise
     \@linelen \hbox to\z@{\lower \ht\@tempboxa\box\@tempboxa\hss}}

\def\@downvector{\@downline\lower \@linelen
      \hbox to \z@{\@linefnt\char'77\hss}}

\def\@getlarrow(#1,#2){\ifnum #2 =\z@ \@tempcnta='33\else
\@tempcnta=#1\relax\multiply\@tempcnta \sixt@@n \advance\@tempcnta
-9 \@tempcntb=#2\relax\multiply\@tempcntb \tw@
\ifnum \@tempcntb >0 \advance\@tempcnta \@tempcntb\relax
\else\advance\@tempcnta -\@tempcntb\advance\@tempcnta 64
\fi\fi\char\@tempcnta}

\def\@getrarrow(#1,#2){\@tempcntb=#2\relax
\ifnum\@tempcntb < 0 \@tempcntb=-\@tempcntb\relax\fi
\ifcase \@tempcntb\relax \@tempcnta='55 \or
\ifnum #1<3 \@tempcnta=#1\relax\multiply\@tempcnta
24 \advance\@tempcnta -6 \else \ifnum #1=3 \@tempcnta=49
\else\@tempcnta=58 \fi\fi\or
\ifnum #1<3 \@tempcnta=#1\relax\multiply\@tempcnta
24 \advance\@tempcnta -3 \else \@tempcnta=51\fi\or
\@tempcnta=#1\relax\multiply\@tempcnta
\sixt@@n \advance\@tempcnta -\tw@ \else
\@tempcnta=#1\relax\multiply\@tempcnta
\sixt@@n \advance\@tempcnta 7 \fi\ifnum #2<0 \advance\@tempcnta 64 \fi
\char\@tempcnta}
%%%%%%%%%%%%%%%%%%%%%%%%%%%%%%%%%%%%%%%%%%%%%%%%%%%%%%%%%%%%%

\def\Diagonale(#1,#2)#3{\unskip\leavevmode
  \xcoord#1\relax \ycoord#2\relax
      \raise\ycoord \Einheit\hbox to0pt{\hskip\xcoord \Einheit
         \Line@(1,1){#3}\hss}}
\def\AntiDiagonale(#1,#2)#3{\unskip\leavevmode
  \xcoord#1\relax \ycoord#2\relax %\advance\xcoord by -0.05\relax
      \raise\ycoord \Einheit\hbox to0pt{\hskip\xcoord \Einheit
         \Line@(1,-1){#3}\hss}}
\def\Pfad(#1,#2),#3\endPfad{\unskip\leavevmode
  \xcoord#1 \ycoord#2 \thicklines\ZeichnePfad#3\endPfad\thinlines}
\def\ZeichnePfad#1{\ifx#1\endPfad\let\next\relax
  \else\let\next\ZeichnePfad
    \ifnum#1=1
      \raise\ycoord \Einheit\hbox to0pt{\hskip\xcoord \Einheit
         \vrule height\Pfadd@cke width1 \Einheit depth\Pfadd@cke\hss}%
      \advance\xcoord by 1
    \else\ifnum#1=2
      \raise\ycoord \Einheit\hbox to0pt{\hskip\xcoord \Einheit
        \hbox{\hskip-\PfadD@cke\vrule height1 \Einheit width\PfadD@cke depth0pt}\hss}%
      \advance\ycoord by 1
    \else\ifnum#1=3
      \raise\ycoord \Einheit\hbox to0pt{\hskip\xcoord \Einheit
         \Line@(1,1){1}\hss}
      \advance\xcoord by 1
      \advance\ycoord by 1
    \else\ifnum#1=4
      \raise\ycoord \Einheit\hbox to0pt{\hskip\xcoord \Einheit
         \Line@(1,-1){1}\hss}
      \advance\xcoord by 1
      \advance\ycoord by -1
    \fi\fi\fi\fi
  \fi\next}
\def\hSSchritt{\leavevmode\raise-.4pt\hbox to0pt{\hss.\hss}\hskip.2\Einheit
  \raise-.4pt\hbox to0pt{\hss.\hss}\hskip.2\Einheit
  \raise-.4pt\hbox to0pt{\hss.\hss}\hskip.2\Einheit
  \raise-.4pt\hbox to0pt{\hss.\hss}\hskip.2\Einheit
  \raise-.4pt\hbox to0pt{\hss.\hss}\hskip.2\Einheit}
\def\vSSchritt{\vbox{\baselineskip.2\Einheit\lineskiplimit0pt
\hbox{.}\hbox{.}\hbox{.}\hbox{.}\hbox{.}}}
\def\DSSchritt{\leavevmode\raise-.4pt\hbox to0pt{%
  \hbox to0pt{\hss.\hss}\hskip.2\Einheit
  \raise.2\Einheit\hbox to0pt{\hss.\hss}\hskip.2\Einheit
  \raise.4\Einheit\hbox to0pt{\hss.\hss}\hskip.2\Einheit
  \raise.6\Einheit\hbox to0pt{\hss.\hss}\hskip.2\Einheit
  \raise.8\Einheit\hbox to0pt{\hss.\hss}\hss}}
\def\dSSchritt{\leavevmode\raise-.4pt\hbox to0pt{%
  \hbox to0pt{\hss.\hss}\hskip.2\Einheit
  \raise-.2\Einheit\hbox to0pt{\hss.\hss}\hskip.2\Einheit
  \raise-.4\Einheit\hbox to0pt{\hss.\hss}\hskip.2\Einheit
  \raise-.6\Einheit\hbox to0pt{\hss.\hss}\hskip.2\Einheit
  \raise-.8\Einheit\hbox to0pt{\hss.\hss}\hss}}
\def\SPfad(#1,#2),#3\endSPfad{\unskip\leavevmode
  \xcoord#1 \ycoord#2 \ZeichneSPfad#3\endSPfad}
\def\ZeichneSPfad#1{\ifx#1\endSPfad\let\next\relax
  \else\let\next\ZeichneSPfad
    \ifnum#1=1
      \raise\ycoord \Einheit\hbox to0pt{\hskip\xcoord \Einheit
         \hSSchritt\hss}%
      \advance\xcoord by 1
    \else\ifnum#1=2
      \raise\ycoord \Einheit\hbox to0pt{\hskip\xcoord \Einheit
        \hbox{\hskip-2pt \vSSchritt}\hss}%
      \advance\ycoord by 1
    \else\ifnum#1=3
      \raise\ycoord \Einheit\hbox to0pt{\hskip\xcoord \Einheit
         \DSSchritt\hss}
      \advance\xcoord by 1
      \advance\ycoord by 1
    \else\ifnum#1=4
      \raise\ycoord \Einheit\hbox to0pt{\hskip\xcoord \Einheit
         \dSSchritt\hss}
      \advance\xcoord by 1
      \advance\ycoord by -1
    \fi\fi\fi\fi
  \fi\next}
\def\Koordinatenachsen(#1,#2){\unskip
 \hbox to0pt{\hskip-.5pt\vrule height#2 \Einheit width.5pt depth1 \Einheit}%
 \hbox to0pt{\hskip-1 \Einheit \xcoord#1 \advance\xcoord by1
    \vrule height0.25pt width\xcoord \Einheit depth0.25pt\hss}}
\def\Koordinatenachsen(#1,#2)(#3,#4){\unskip
 \hbox to0pt{\hskip-.5pt \ycoord-#4 \advance\ycoord by1
    \vrule height#2 \Einheit width.5pt depth\ycoord \Einheit}%
 \hbox to0pt{\hskip-1 \Einheit \hskip#3\Einheit
    \xcoord#1 \advance\xcoord by1 \advance\xcoord by-#3
    \vrule height0.25pt width\xcoord \Einheit depth0.25pt\hss}}
\def\Gitter(#1,#2){\unskip \xcoord0 \ycoord0 \leavevmode
  \LOOP\ifnum\ycoord<#2
    \loop\ifnum\xcoord<#1
      \raise\ycoord \Einheit\hbox to0pt{\hskip\xcoord \Einheit\Punkt\hss}%
      \advance\xcoord by1
    \repeat
    \xcoord0
    \advance\ycoord by1
  \REPEAT}
\def\Gitter(#1,#2)(#3,#4){\unskip \xcoord#3 \ycoord#4 \leavevmode
  \LOOP\ifnum\ycoord<#2
    \loop\ifnum\xcoord<#1
      \raise\ycoord \Einheit\hbox to0pt{\hskip\xcoord \Einheit\Punkt\hss}%
      \advance\xcoord by1
    \repeat
    \xcoord#3
    \advance\ycoord by1
  \REPEAT}
\def\Label#1#2(#3,#4){\unskip \xdim#3 \Einheit \ydim#4 \Einheit
  \def\lo{\advance\xdim by-.5 \Einheit \advance\ydim by.5 \Einheit}%
  \def\llo{\advance\xdim by-.25cm \advance\ydim by.5 \Einheit}%
  \def\loo{\advance\xdim by-.5 \Einheit \advance\ydim by.25cm}%
  \def\o{\advance\ydim by.25cm}%
  \def\ro{\advance\xdim by.5 \Einheit \advance\ydim by.5 \Einheit}%
  \def\rro{\advance\xdim by.25cm \advance\ydim by.5 \Einheit}%
  \def\roo{\advance\xdim by.5 \Einheit \advance\ydim by.25cm}%
  \def\l{\advance\xdim by-.30cm}%
  \def\r{\advance\xdim by.30cm}%
  \def\lu{\advance\xdim by-.5 \Einheit \advance\ydim by-.6 \Einheit}%
  \def\llu{\advance\xdim by-.25cm \advance\ydim by-.6 \Einheit}%
  \def\luu{\advance\xdim by-.5 \Einheit \advance\ydim by-.30cm}%
  \def\u{\advance\ydim by-.30cm}%
  \def\ru{\advance\xdim by.5 \Einheit \advance\ydim by-.6 \Einheit}%
  \def\rru{\advance\xdim by.25cm \advance\ydim by-.6 \Einheit}%
  \def\ruu{\advance\xdim by.5 \Einheit \advance\ydim by-.30cm}%
  #1\raise\ydim\hbox to0pt{\hskip\xdim
     \vbox to0pt{\vss\hbox to0pt{\hss$#2$\hss}\vss}\hss}%
}

\font@\twelverm=cmr10 scaled\magstep1
\font@\twelveit=cmti10 scaled\magstep1
\font@\twelvebf=cmbx10 scaled\magstep1
\font@\twelvei=cmmi10 scaled\magstep1
\font@\twelvesy=cmsy10 scaled\magstep1
\font@\twelveex=cmex10 scaled\magstep1

\newtoks\twelvepoint@
\def\twelvepoint{\normalbaselineskip15\p@
 \abovedisplayskip15\p@ plus3.6\p@ minus10.8\p@
 \belowdisplayskip\abovedisplayskip
 \abovedisplayshortskip\z@ plus3.6\p@
 \belowdisplayshortskip8.4\p@ plus3.6\p@ minus4.8\p@
 \textonlyfont@\rm\twelverm \textonlyfont@\it\twelveit
 \textonlyfont@\sl\twelvesl \textonlyfont@\bf\twelvebf
 \textonlyfont@\smc\twelvesmc \textonlyfont@\tt\twelvett
%Ergnzung des fetten Small-Capitals-Fonts:
%
 \ifsyntax@ \def\big##1{{\hbox{$\left##1\right.$}}}%
  \let\Big\big \let\bigg\big \let\Bigg\big
 \else
  \textfont\z@=\twelverm  \scriptfont\z@=\tenrm  \scriptscriptfont\z@=\sevenrm
  \textfont\@ne=\twelvei  \scriptfont\@ne=\teni  \scriptscriptfont\@ne=\seveni
  \textfont\tw@=\twelvesy \scriptfont\tw@=\tensy \scriptscriptfont\tw@=\sevensy
  \textfont\thr@@=\twelveex \scriptfont\thr@@=\tenex
        \scriptscriptfont\thr@@=\tenex
  \textfont\itfam=\twelveit \scriptfont\itfam=\tenit
        \scriptscriptfont\itfam=\tenit
  \textfont\bffam=\twelvebf \scriptfont\bffam=\tenbf
        \scriptscriptfont\bffam=\sevenbf
  \setbox\strutbox\hbox{\vrule height10.2\p@ depth4.2\p@ width\z@}%
  \setbox\strutbox@\hbox{\lower.6\normallineskiplimit\vbox{%
        \kern-\normallineskiplimit\copy\strutbox}}%
 \setbox\z@\vbox{\hbox{$($}\kern\z@}\bigsize@=1.4\ht\z@
 \fi
 \normalbaselines\rm\ex@.2326ex\jot3.6\ex@\the\twelvepoint@}

\font@\fourteenrm=cmr10 scaled\magstep2
\font@\fourteenit=cmti10 scaled\magstep2
\font@\fourteensl=cmsl10 scaled\magstep2
\font@\fourteensmc=cmcsc10 scaled\magstep2
\font@\fourteentt=cmtt10 scaled\magstep2
\font@\fourteenbf=cmbx10 scaled\magstep2
\font@\fourteeni=cmmi10 scaled\magstep2
\font@\fourteensy=cmsy10 scaled\magstep2
\font@\fourteenex=cmex10 scaled\magstep2
\font@\fourteenmsa=msam10 scaled\magstep2
\font@\fourteeneufm=eufm10 scaled\magstep2
\font@\fourteenmsb=msbm10 scaled\magstep2
\newtoks\fourteenpoint@
\def\fourteenpoint{\normalbaselineskip15\p@
 \abovedisplayskip18\p@ plus4.3\p@ minus12.9\p@
 \belowdisplayskip\abovedisplayskip
 \abovedisplayshortskip\z@ plus4.3\p@
 \belowdisplayshortskip10.1\p@ plus4.3\p@ minus5.8\p@
 \textonlyfont@\rm\fourteenrm \textonlyfont@\it\fourteenit
 \textonlyfont@\sl\fourteensl \textonlyfont@\bf\fourteenbf
 \textonlyfont@\smc\fourteensmc \textonlyfont@\tt\fourteentt
%Erg^=C4nzung des fetten Small-Capitals-Fonts:
%
 \ifsyntax@ \def\big##1{{\hbox{$\left##1\right.$}}}%
  \let\Big\big \let\bigg\big \let\Bigg\big
 \else
  \textfont\z@=\fourteenrm  \scriptfont\z@=\twelverm  \scriptscriptfont\z@=\tenrm
  \textfont\@ne=\fourteeni  \scriptfont\@ne=\twelvei  \scriptscriptfont\@ne=\teni
  \textfont\tw@=\fourteensy \scriptfont\tw@=\twelvesy \scriptscriptfont\tw@=\tensy
  \textfont\thr@@=\fourteenex \scriptfont\thr@@=\twelveex
        \scriptscriptfont\thr@@=\twelveex
  \textfont\itfam=\fourteenit \scriptfont\itfam=\twelveit
        \scriptscriptfont\itfam=\twelveit
  \textfont\bffam=\fourteenbf \scriptfont\bffam=\twelvebf
        \scriptscriptfont\bffam=\tenbf
  \setbox\strutbox\hbox{\vrule height12.2\p@ depth5\p@ width\z@}%
  \setbox\strutbox@\hbox{\lower.72\normallineskiplimit\vbox{%
        \kern-\normallineskiplimit\copy\strutbox}}%
 \setbox\z@\vbox{\hbox{$($}\kern\z@}\bigsize@=1.7\ht\z@
 \fi
 \normalbaselines\rm\ex@.2326ex\jot4.3\ex@\the\fourteenpoint@}

\catcode`\@=13
\def\[{\left[}
\def\]{\right]}
\define\twoline#1#2{\line{\hfill{\smc #1}\hfill{\smc #2}\hfill}}

\def\mypic#1{\epsffile{#1}}

\topmatter
\title The interaction of a gap with a free boundary in a 
two dimensional dimer system
\endtitle
\author M.~Ciucu$^\dagger$ and C.~Krattenthaler$^\ddagger$
\endauthor
\affil Department of Mathematics,
Indiana University,\\
Bloomington, IN 47405-5701, USA\\\vskip6pt
Fakult\"at f\"ur Mathematik der Universit\"at Wien,\\
Nordbergstra{\ss}e 15, A-1090 Wien, Austria.\\
WWW: \tt http://www.mat.univie.ac.at/\~{}kratt
\endaffil
\address Department of Mathematics, Indiana University, Bloomington, IN 47405-5701,
USA
\endaddress
\address Fakult\"at f\"ur Mathematik der Universit\"at Wien,
Nordbergstrasze 15, A-1090 Wien, Austria.
\endaddress
%\dedicatory \enddedicatory
%\date \enddate
\thanks{$^\dagger$Research partially supported by NSF grant DMS-0500616.}\endthanks
\thanks{$^\ddagger$Research partially supported by the Austrian
Science Foundation FWF, grants Z130-N13 and S9607-N13,
the latter in the framework of the National Research Network
``Analytic Combinatorics and Probabilistic Number Theory."}\endthanks
%\subjclass[2000] Primary 05A15;
% Secondary 05A16 05A17 05A19 05B45 33C20 52C20
%\endsubjclass
%\keywords rhombus tilings, lozenge tilings, plane partitions,
%nonintersecting lattice paths, determinant evaluations\endkeywords
\abstract
Let $\ell$ be a fixed vertical lattice line of the unit triangular
lattice in the plane, and let $\Cal H$ be the half  
plane to the left of $\ell$. We consider lozenge tilings of $\Cal H$
that have a triangular gap of side-length two 
and in which $\ell$ is a free boundary --- i.e., tiles are allowed to
protrude out half-way across $\ell$. We prove  
that the correlation function of this gap near the free boundary has
asymptotics $\frac{1}{4\pi r}$, $r\to\infty$, 
where $r$ is the distance from the gap to the free boundary. This
parallels the electrostatic phenomenon by which the 
field of an electric charge near a conductor can be obtained by the
method of images. 
\endabstract
\endtopmatter
\document

\leftheadtext{M. Ciucu and C. Krattenthaler}

\rightheadtext{The interaction of a gap with a free boundary in a
dimer system}

\head 1. Introduction\endhead
The study of the interaction of gaps in dimer coverings was introduced
in the literature by Fisher 
and Stephenson \cite{\FS}. This pioneering work contains three
different types of gap interaction  
in dimer systems on the square lattice:
$(i)$ interaction of two dimer-gaps (equivalently, interaction of two
fixed dimers required to be  
contained in the dimer coverings); 
$(ii)$ interaction of two non-dimer-gaps (specifically, two monomers), 
and $(iii)$ the interaction of a dimer-gap
with a constrained boundary (edge or corner). 

The first of these types of interactions was later 
generalized by Kenyon \cite{\Kone} to an arbitrary number of
dimer-gaps on the square and  
hexagonal lattices, and recently by Kenyon, Okounkov and Sheffield
\cite{\KOS} to general planar bipartite lattices. 
Interactions of the second type were studied by the first author of the present
paper in \cite{\CiucAK}\cite{\sc}\cite{\CiucAM}\cite{\ef}\cite{\ov},
where close analogies to two dimensional electrostatics were established. 

Two instances of interaction of non-dimer-gaps with constrained
boundaries can be found in \cite{\Kmonbd, Section~7.5} (interaction
of a monomer with a constrained straight line boundary on 
the square lattice), and respectively \cite{\sc, Theorem~2.2}
(interaction of a family of triangular gaps with a constrained
straight line boundary on the hexagonal lattice). 

In this paper we determine the interaction of a triangular gap with a
{\it free straight line boundary} (i.e., dimers are allowed to
protrude out across it) on the hexagonal lattice. 
This type of interaction seems not to have been treated before
in the literature. 
(We are aware of one other paper, namely \cite{\FrReAA}, 
addressing the asymptotic behavior of lozenge tilings
under the presence of a free boundary, but the regions considered
there contain no gaps.) We find that the gap is attracted to the free
boundary in precise analogy to the (two dimensional) electrostatic
phenomenon in which an electric charge is attracted by a  
straight line conductor when placed near it.

This develops further the analogy between dimer systems with gaps and
electrostatics that the first author has described in
\cite{\sc}\cite{\ec}\cite{\ef}\cite{\ov}. More generally, our result
shows that in any physical system that can be modeled by dimer
coverings, a gap will tend to be attracted to an interface
corresponding to a free boundary. This effect, purely entropic in
origin, is reminiscent of the Cheerios effect by which an air bubble
at the surface of a liquid in a container is attracted to the walls
\cite{\VM} (note that the Cheerios effect is {\it not} entropic in
origin). 

\head 2. Set-up and results\endhead
There seem to be no methods in the literature for finding the
interaction of a gap ``in a sea of dimers'' with a free
boundary. However, as V.~I.~Arnold said, ``mathematics is a part of
physics where experiments are cheap.'' We now design such an
experiment in order to determine the interaction of a gap in a dimer
system on the hexagonal lattice with a free boundary.

Consider the tiling of the plane by unit triangles, drawn so that one
family of lattice lines is vertical. Clearly, the hexagonal lattice
can be viewed as the graph whose vertices are the unit triangles, and
whose edges connect precisely those unit triangles that share an
edge. Dimers on the hexagonal lattice are then (unit) lozenges (i.e.,
unit rhombi) consisting of pairs of adjacent unit triangles.

The free boundary we choose is a lattice line $\ell$ --- say vertical
--- on the triangular lattice, to the left of which
the plane is covered completely and without overlapping by lozenges,
except for a gap $\triangleleft_2$ in the shape  
of a triangle of side-length 2, pointing to the left; the lozenges are
allowed to protrude halfway across the free  
boundary, to its right (Figure~\FF\ pictures a portion of such a
tiling; the dotted lines should be ignored at this point).

We define the {\it correlation function} (or simply correlation) of
the hole $\triangleleft_2$ with the free  
boundary $\ell$ as follows. Choose a rectangular system of coordinates
in which $\ell$ is the $y$-axis, the origin 
is some lattice point on $\ell$, and the unit is the lattice spacing. 
Let $\triangleleft_2(k)$ be the placement of $\triangleleft_2$ so that
the center $C$ of its right side has  
coordinates $(-k\sqrt{3},0)$ (i.e., $C$ and the origin are the
endpoints of a string of $k$ contiguous horizontal  
lozenges; Figure~\FF\ illustrates $\triangleleft_2(2)$, the origin
being denoted by $O$ there). 
Let $H_{n,x}$ be the lattice hexagon of side-lengths $2n$, $2n$, $2x$,
$2n$, $2n$, $2x$ (in counter-clockwise order, 
starting with the southwestern side) centered
at the origin 
(thus $H_{n,x}$ is vertically symmetric about $\ell$, and its
horizontal symmetry axis cuts $\triangleleft_2(k)$  
into two equal parts; for example, 
the boundary of the region in Figure~\FE\ is $H_{4,4}$). 
Let $F_{n,x}$ be the region obtained from the left half of $H_{n,x}$ by regarding
its boundary along $\ell$ as free (i.e., lozenges in a tiling of
$F_{n,x}$ are allowed to protrude outward across 
$\ell$). Figure~\FFa\ shows the region $F_{3,3}$
together with such a lozenge tiling; 
the origin is labelled by $O$. 

\midinsert
\vbox{
\centerline{\mypic{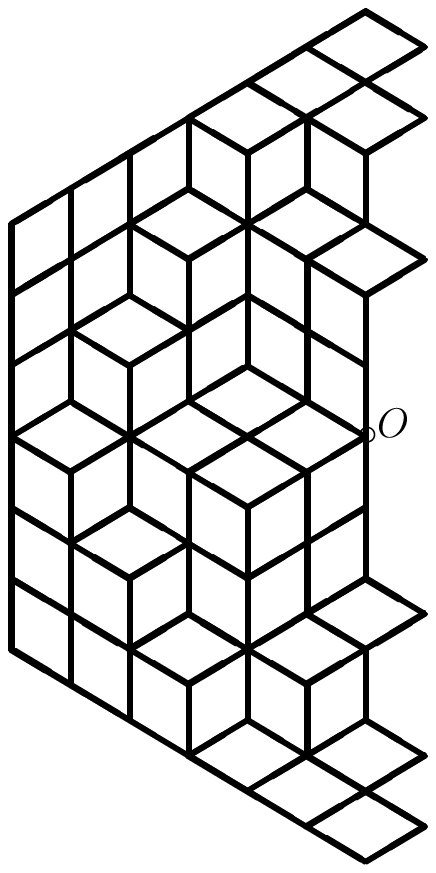}}
\centerline{\eightpoint Figure \FFa}
}
\vskip10pt
\endinsert

Following in the spirit of \cite{\FS} and \cite{\sc}, 
for any fixed integer $k\geq0$ we define the {\it correlation of
$\triangleleft_2(k)$ with the free boundary $\ell$},  
denoted $\omega_f(k)$, by
$$
\omega_f(k):=\lim_{n\to\infty}\frac{\M(F_{n,n}\setminus
\triangleleft_2(k))}{\M(F_{n,n})},\tag\Bmm 
$$
where $\M(R)$ stands for the number of lozenge tilings of the region
$R$ (if $R$ has portions of the boundary that 
are free --- as in our case --- then it is understood that what we
count is tilings in which lozenges are allowed 
to protrude out across the free portions). A tiling of
$F_{4,4}\setminus\triangleleft_2(2)$ of this type
is illustrated in Figure~\FF.

We note that, by \cite{\CLP}, lozenges have maximum entropy statistics
(in the scaling limit) at the center of 
a regular hexagon. According to this, ({\Bmm}) is a natural definition
for the correlation function. An analogous  
definition was used in \cite{\sc}.

In Lemma~{\TD} (with $\xi=1$) 
we obtain an exact expression for $\omega_f(k)$ in terms of an
integral. What affords this is 
an exact formula for $\M(F_{n,x}\setminus \triangleleft_2(k))$, which we present in
Theorem~\TB. We then deduce 
the asymptotics of $\omega_f(k)$ as $k\to\infty$ using Laplace's
method (see Theorem~\TF\ with $\xi=1$). The result is the following. 

\proclaim{Theorem \TAA} As $k\to\infty$, we have
%$(${\rm a}$)$.
$$
\omega_f(k)\sim\frac{1}{4\pi}\frac{1}{\de(\triangleleft_2(k),\ell)},\tag\Bm
$$
where $\de$ is the Euclidean distance.
\endproclaim

\remark{Remark 1}
In fact, our results allow us to determine the correlation of
$\triangleleft_2(k)$ with the free boundary $\ell$ 
in a more general situation, namely when the sides $2n$ and $2x$ of
$F_{n,x}$ grow to infinity so that $x/n$ approaches a positive real
number $\xi$ not necessarily equal to $1$. This leads to the
correlation $\omega_f(k;\xi)$ defined in (\X). It turns out (see
Theorem~\TF) that, for $\xi\ne1$, the behavior of this correlation is
{\it exponential\/} in $k$. This is in contrast to the behavior of
the correlation of lozenges on the (``infinitely large'') torus, in which case 
Kenyon, Okounkov and Sheffield have shown in \cite{\KOS, Sec.~4.4}
that the correlation decays {\it polynomially}. More precisely,
it is shown in \cite{\KOS} that, for dimer models on doubly periodic 
bipartite planar graphs, 
there can only occur three different {\it``phases"} characterized by the
behavior of edge-edge correlations: liquid, gaseous, and frozen.  
Our situation is readily seen to be in the liquid phase. The liquid phase
is shown in \cite{\KOS} to be characterized by a polynomial decay of edge-edge 
correlations, and the arguments there imply that the correlation of holes of 
side 2 also have polynomial asymptotic behavior. It is this fact that is in
contrast with the exponential interaction of Theorem~\TF.
In fact, if one
would define analogous generalized correlations for the lozenge tiling
models in \cite{\CiucAK}\cite{\sc}\cite{\CiucAM}, one would observe the same
phenomenon of exponential behavior of correlation if the correlation
is not ``measured" inside a patch of the lattice region in which all
three types of lozenges occur with equal probability. This hints at
the limitation of the transfer of properties of dimer models on the
torus to dimer models of bounded regions in the plane, which is one of
the main driving forces in \cite{\KOS}. We plan to address this
phenomenon in more detail in forthcoming publications.
%In contrast to this, the correlation of a lozenge
%with a free boundary is expected to decay polynomially in $k$, based
%on the results of Kenyon, Okounkov and Sheffield in \cite{\KOS,
%Sec.~4.4}. (The correlation of the lozenge $L$ with a free boundary $\ell$ is
%expected to be determined by the correlation of $L$ with its mirror
%image $L'$ with respect to $\ell$ defined now by enclosing $L$ and
%$L'$ in the full hexagon $H_{n,x}$. By the results in \cite{\KOS}, the
%latter decays polynomially.) This points to the fact that the
%correlation of holes of nonzero charge\footnote{The {\it charge} of a
%hole is equal to the number of right-pointing unit triangles inside it
%minus the number of the corresponding left-pointing ones.} behave
%radically differently than the correlation of zero charge holes 
%(such as, for instance, lozenge holes). 
\endremark

\medskip
In \cite{\ov} the first author described how a distribution of fixed
holes on the triangular lattice defines in a natural way 
two vector fields. The ${\bold F}$-field is a discrete vector field
defined at the center of each left-pointing unit 
triangle $e$, and equal to the expected orientation of the lozenge
covering $e$ (under the uniform measure on the 
set of tilings). To define the ${\bold T}$-field, one introduces an
extra ``test-hole'' $t$ and measures the relative 
change in the correlation function under small displacements of it, as
the other holes are kept fixed.  
One can prove (details will appear elsewhere) 
that in the scaling limit of the lattice spacing approaching zero, 
this relative change is given by the scalar product of the displacement vector with
a certain vector ${\bold T}(z)$, where $z$ is the point to which the
test hole $t$ shrinks when the lattice spacing 
approaches zero. This defines the second field.

When these fields are generated by lozenge tilings that
cover the entire plane with the exception of a finite collection of
fixed-size holes (the case treated in \cite{\ov} and \cite{\ef}), both
the ${\bold T}$-field and the scaling limit of the ${\bold F}$-field
turn out to be equal, up to a constant multiple, to the electrostatic
field of the two dimensional physical system obtained by viewing the
holes as electrical charges. 

But what if we do not tile the entire plane, but only the half-plane to
the left of the free boundary $\ell$, and we have no holes?

The above definitions for the ${\bold F}$-field and ${\bold T}$-field
would still work, provided $(i)$ the  
scaling limit of the discrete field defining ${\bold F}$ exists, and
$(ii)$ the scaling limit of the relative 
changes in the correlation function under small displacements of a
test hole exists and is given by taking  
scalar products of the displacement vector with the vectors of a certain field.

Our exact determination of $\omega_f(k)$ (see Lemma~{\TD}) allows us
to verify $(ii)$ for displacements along the horizontal direction. 
%draw an additional conclusion from our ``experiment.''
$\triangleleft_2(k)$ plays now the role of a test charge. We obtain
the following result. 
%Note that moving the hole from $\triangleleft_2(k)$ to
%$\triangleleft_2(k+\alpha)$ amounts to displacing it by a distance of $\alpha\sqrt{3}$.

\proclaim{Theorem \TAB} We have
$$
%\frac{1}{2\alpha}
%\left(
\frac{\omega_f(k+1)}{\omega_f(k)}-1
%\right)
\sim
%-\frac{1}{2\de(\triangleleft_2(k),\ell)},
-\frac{1}{k},
\ \ \ k\to\infty.\tag\TABzero
$$

\endproclaim

\remark{Remark 2}
By symmetry, displacements of $\triangleleft_2(k)$
parallel to $\ell$ leave $\omega_f$ unchanged, 
so the relative change in $\omega_f$ corresponding to such
displacements is zero. Thus, provided the  
field ${\bold T}$ exists, it follows from Theorem~{\TAB} that its
value at $z$ is 
$$
{\bold T}(z)=-\frac{\bold i}{2\de(z,\ell)},\tag\TABone
$$
where ${\bold i}$ is the unit vector in the positive direction of the
$x$-axis (the 2 at the denominator comes 
from the fact that ${\bold T}$ arises from the expression on the left
hand side of (\TABzero) divided by the product 
of the displacement, $\sqrt{3}$ in this case, and the ``charge''
of the hole $\triangleleft_2(k)$, which is 2; see \cite{\ov} for
details). Note that by \cite{\ov} we 
would obtain (up to a multiplicative constant of 2) the same field
${\bold T}$ at $z$ if we look at tilings of the  
entire plane, with the mirror image of
our test-hole $\triangleleft_2(k)$ being a fixed hole.
This is analogous to the phenomenon
in electrostatics by which the field created by an electric charge
placed near a conductor can be obtained  
by the method of images (see e.g. \cite{\Feytwo, Chapter~6}).
\endremark

\medskip
The ${\bold F}$-field could be determined by an ``experiment''
analogous to the one we described at the beginning of this
section. What one needs now is the number of lozenge tilings of
$F_{n,x}\setminus L(k)$, where $L(k)$ is the horizontal lozenge
contained in $\triangleleft_2(k)$. This turns out to be given by a
formula very similar to ({\ABa}), namely by 
$$
\multline
\M(F_{n,x}\setminus L(k))=
{\binom {n+k}{2k+1}}^{\!2}
\frac {(n-k-1)!\,(x+n-k)_{2k+1}} {(n-k)_{2k+1}}
\prod _{s=1} ^{n}\frac {(2x+2s)_{4n-4s+1}} {(2s)_{4n-4s+1}}
\\
\times
\sum _{i=0} ^{n-k-1}
\frac{1}{{\binom {n+k-i}{2k+1}}^{\!2}}
\frac {    (\tfrac 12)_{i}}
  {i!\,(n-k-i-1)!^2\,(n+k-i+1)_{n-k}\,(n+k-i+1)_i\,(2n-i+\frac {1} {2})_i}\\
\cdot
\Big((x)_{i}\,(x+i+1)_{n-k-i-1}\,(x+n+k+1)_{n-k}\\
-
          (x)_{n-k}\,(x+n+k+1)_{n-k-i-1}\,(x+2n-i+1)_{i}\Big).
\endmultline\tag\Bp
$$

\flushpar
In the same way as we derive Theorem~\TAA\ from Theorem~\TB\ via
Lemmas~\TC--\TE, 
Laplace's method can be used to deduce from (\Bp) the following result.

\proclaim{Corollary \TAC} Let $e(k)$ be the leftmost left-pointing
unit triangle of $\triangleleft_2(k)$. Then 
$$
{\bold F}(e(k))\sim o\left(\frac{1}{\de(e(k),\ell)}\right),\ \ \ k\to\infty.\tag\Bpp
$$
\endproclaim

\remark{Remark 3}
Equations ({\TABone}) and ({\Bpp}) imply that, in sharp contrast to the case
of lozenge tilings of the plane with a finite number of fixed size
holes, for the half-plane with free boundary the 
fields ${\bold T}$ and ${\bold F}$ have radically different behavior:
while the former behaves as the electrostatic 
field near a conductor, the latter is zero in the scaling limit. 
\endremark

Our approach to proving Theorems~\TAA\ and \TAB\ consists of solving
first the counting problem exactly, see Theorem~\TB. This result
generalizes Andrews' theorem \cite{\AndrAK} (which proved MacMahon's
conjecture on symmetric plane partitions) in the case $q=1$.
Its proof is given in Sections~4 and 5, with some auxiliary results
proved separately in Section~6. It is based on the
``exhaustion/identification of factors" method described in
\cite{\KratBN, Sec.~2.4}. Finally, in Section~7, we perform the
asymptotic calculations needed to derive Theorems~\TAA\ and \TAB\ from
the exact counting results.

\head 3. An exact tiling enumeration formula\endhead

Tilings of the region $F_{n,x}$ are clearly equivalent to tilings of
the hexagon $H_{n,x}$ that are invariant under  
reflection across its symmetry axis $\ell$. Counting such tilings was
a problem considered (in the equivalent form of 
symmetric plane partitions)
by MacMahon in the early twentieth century (see \cite{\MacMAA, p.~270}). 
MacMahon conjectured that the number of vertically symmetric lozenge tilings of 
a hexagon with side lengths $2n,2n,2x,2n,2n,2x$ is equal to
$$\frac {\(x+\frac {1} {2}\)_{2n}} {\(\frac {1} {2}\)_{2n}}
{\prod _{s=1} ^{n}}\frac {(2x+2s)_{4n-4s+1}} 
{(2s)_{4n-4s+1}},\tag\SPP
$$
where $(\alpha)_m$ is the Pochhammer symbol, defined by
$(\alpha)_m:=\alpha(\alpha+1)\cdots(\alpha+m-1)$ for 
$m\ge1$, and $(\alpha)_0:=1$.
This was first proved
by Andrews \cite{\AndrAK}. Other proofs,
and refinements, were later found by e.g.\ 
Gordon \cite{\GordAC},
Macdonald \cite{\MacdAC, pp.~83--85}, Proctor \cite{\ProcAD,
Prop.~7.3}, Fischer \cite{\FiscAH}, and  
the second author of the present paper \cite{\KratAP, Theorem~13}.

Our ``experiment'' --- counting
$\M(F_{n,x}\setminus\triangleleft_2(k))$ --- is by the same token
equivalent to 
counting vertically symmetric lozenge tilings of $H_{n,x}$ with two
missing triangles (compare Figures~\FE\ and \FF). This is in fact a
generalization of MacMahon's symmetric plane partitions problem (see Remark~4).

\medskip
The key result that allows deducing Theorems~{\TAA} and {\TAB} is the following.

\proclaim{Theorem \TB}
For all positive integers $n,x$ and nonnegative integers $k\leq n-1$, we have
$$\multline 
\M(F_{n,x}\setminus\triangleleft_2(k))=
\binom {4k+1}{2k}\,
\frac {(n+k)!} {(x+n-k)_{2k+1}}
{\prod _{s=1} ^{n}}\frac {(2x+2s)_{4n-4s+1}} 
{(2s)_{4n-4s+1}}\\
\times
\sum _{i=0} ^{n-k-1}
\frac {    (\tfrac 12)_{i}}
  {i!\,(n-k-i-1)!^2\,(n+k-i+1)_{n-k}\,(n+k-i+1)_i\,(2n-i+\frac {1} {2})_i}\\
\cdot
\Big((x)_{i}\,(x+i+1)_{n-k-i-1}\,(x+n+k+1)_{n-k}\\
-
          (x)_{n-k}\,(x+n+k+1)_{n-k-i-1}\,(x+2n-i+1)_{i}\Big).
%
%\binom {4k+1}{2k}\,
%\frac {1} {(n+k+\tfrac 32)_{n-k-1}\,
%     (2k+2)_{n-k-1}}
%\underset s\ne n-k\to{\prod _{s=1} ^{n}}\frac {(x+s)_{2n-2s+1}} 
%{(s)_{2n-2s+1}}
%{\prod _{s=1} ^{n}}\frac {(x+s+\tfrac {1} {2})_{2n-2s}}
%{(s+\tfrac {1} {2})_{2n-2s}}\\
%\times
%\sum _{i=0} ^{n-k-1}
%\frac {     (n-k-i)_{i}\,(\tfrac 12)_{i}\,(n+k+\frac {\scriptstyle3} {\scriptstyle2})_{n-k-i-1}\,
%     (2k+2)_{n-k-i-1}}
%  {i!\,(n-k-i-1)!\,(n+k+1-i)_{n-k}}\hskip4cm\\
%\cdot
%\Big((x)_{i}\,(x+i+1)_{n-k-i-1}\,(x+n+k+1)_{n-k}\\
%-
%          (x)_{n-k}\,(x+n+k+1)_{n-k-i-1}\,(x+2n-i+1)_{i}\Big).
\endmultline\tag\ABa$$

\endproclaim

\remark{Remark 4}
Replacing $x$ by $x-1$, $n$ by $n+1$, and $k$ by $n$, one can see that
the above formula specializes to MacMahon's formula (\SPP). More
precisely, because of forced lozenges (see Figure~{\FC}), the
enumeration problem in the statement of Theorem~\TB\ reduces to the
problem of enumerating vertically symmetric lozenge tilings of a
hexagon with side lengths $2n,2n,2x,2n,2n,2x$.
\endremark

\topinsert
\vskip10pt
\vbox{\noindent
\centertexdraw{
\drawdim truecm \setunitscale.8
\linewd.15
\RhombusC \RhombusC 
\move(0 1)
\RhombusC \RhombusC 
\move(0 2)
\RhombusC \RhombusC 
\move(0 3)
\RhombusC \RhombusC 
\move(0 7)
\RhombusB \RhombusB \RhombusB \RhombusB 
\move(0.866025 7.5)
\RhombusB \RhombusB \RhombusB \RhombusB 
\linewd.02
\move(1.73205 -1)
\rdreieck \rhombus \rhombus \rhombus 
\move(1.73205 0)
\rdreieck \rhombus \rhombus \rhombus 
\move(1.73205 1)
\rdreieck \rhombus \rhombus \rhombus 
\move(1.73205 1)
\rhombus \rhombus \rhombus 
\move(1.73205 3)
\rdreieck \rhombus \rhombus \rhombus 
\move(1.73205 4)
\rdreieck \rhombus \rhombus \rhombus 
\move(1.73205 4)
\rhombus \rhombus \rhombus 
\move(1.73205 5)
\rhombus \rhombus \rhombus 
\move(1.73205 6)
\rhombus \rhombus \rhombus 
\move(1.73205 7)
\rhombus \rhombus \rhombus 
\move(1.73205 8)
\rhombus \rhombus \rhombus 
\move(2.59808 8.5)
\rhombus \rhombus
\move (3.4641 9)
\rhombus
\move(4.33013 -3.5)
\huSchritt
\move(4.33013 9.5)
\hoSchritt
\move(5.19615 -3)
\vdSchritt
\move(5.19615 -2)
\vdSchritt
\move(5.19615 -1)
\vdSchritt
\move(5.19615 0)
\vdSchritt
\move(5.19615 1)
\vdSchritt
\move(5.19615 2)
\vdSchritt
\move(5.19615 3)
\vdSchritt
\move(5.19615 4)
\vdSchritt
\move(5.19615 5)
\vdSchritt
\move(5.19615 6)
\vdSchritt
\move(5.19615 7)
\vdSchritt
\move(5.19615 8)
\vdSchritt
\move(5.19615 9)
\vdSchritt
\move(5.19615 10)
\vdSchritt
}
\centerline{\eightpoint Forced lozenges when the hole touches the left border}
\vskip8pt
\centerline{\eightpoint Figure \FC}
}
\vskip10pt
\endinsert

The proof of Theorem~\TB\ is given in the next two sections.
In Section~4, we show that $\M(F_{n,x}\setminus\triangleleft_2(k))$
can be expressed in terms of a certain Pfaffian. This Pfaffian is then
evaluated in Section~5.

\head 4. Lozenge tilings and nonintersecting lattice paths\endhead

\topinsert
\vskip10pt
\vbox{\noindent
\centertexdraw{
\drawdim truecm \setunitscale.7
\linewd.15
\move(0 4)
\RhombusB \RhombusB \RhombusB 
\RhombusA \RhombusB \RhombusA \RhombusB \RhombusA \RhombusA \RhombusB 
\RhombusA \RhombusB \RhombusB \RhombusA \RhombusA \RhombusA
\move(0.866025 4.5)
\RhombusB \RhombusB \RhombusB \RhombusA \RhombusB
\move(1.73205 5)
\RhombusB \RhombusB \RhombusB \RhombusA \RhombusB
\move(2.59808 5.5)
\RhombusB \RhombusB \RhombusA \RhombusB \RhombusB \RhombusA \RhombusB \RhombusB
\RhombusA \RhombusB \RhombusA \RhombusB \RhombusA 
\move(3.4641 6)
\RhombusB \RhombusB \RhombusA \RhombusB \RhombusB \RhombusA \RhombusA \RhombusB
\RhombusB \RhombusB \RhombusA 
\move(4.33013 6.5)
\RhombusB \RhombusB \RhombusA \RhombusB \RhombusA \RhombusB
\RhombusA
\move(5.19615 7)
\RhombusB \RhombusA \RhombusB \RhombusB \RhombusA
\move(6.06218 7.5)
\RhombusB \RhombusA
\move(0 1)
\RhombusC 
\move(0 0)
\RhombusC \RhombusC 
\move(0 -1)
\RhombusC \RhombusC \RhombusC \RhombusC 
\move(0 -2)
\RhombusC \RhombusC \RhombusC \RhombusC \RhombusC 
\move(0 -3)
\RhombusC \RhombusC \RhombusC \RhombusC \RhombusC 
\move(3.4641 1)
\RhombusC 
\move(3.4641 0)
\RhombusC 
\move(2.59808 -.5)
\RhombusC \RhombusC \RhombusC 
\move(4.33013 -2.5)
\RhombusC \RhombusC 
\move(5.19615 -4)
\RhombusC \RhombusC 
\move(5.19615 -5)
\RhombusC \RhombusC 
\move(5.19615 1)
\RhombusC 
%symmetry
\move(6.9282 8)
\RhombusC \RhombusC \RhombusC \RhombusC 
\RhombusC \RhombusC \RhombusC \RhombusC 
\move(6.9282 6)
\RhombusA \RhombusB \RhombusA \RhombusA \RhombusA 
\RhombusB \RhombusA \RhombusA \RhombusA 
\RhombusB \RhombusB \RhombusB \RhombusB \RhombusB 
\move(6.9282 5)
\RhombusC \RhombusC \RhombusC \RhombusC 
\move(6.9282 3)
\RhombusA \RhombusB \RhombusA \RhombusA \RhombusB  \RhombusB \RhombusB 
\RhombusA \RhombusB  \RhombusB \RhombusB 
\move(6.9282 1)
\RhombusA \RhombusB \RhombusB \RhombusA \RhombusB  \RhombusB  
\RhombusA \RhombusB  \RhombusB 
\move(6.9282 1)
\RhombusC 
\move(6.9282 0)
\RhombusC 
\move(6.9282 -1)
\RhombusC \RhombusC 
\move(6.9282 -3)
\RhombusA \RhombusB \RhombusB \RhombusB \RhombusA
\move(6.9282 -5)
\RhombusB \RhombusB \RhombusA
\move(8.66025 -4)
\RhombusC
\move(8.66025 -5)
\RhombusC
\move(8.66025 0)
\RhombusC
\move(8.66025 1)
\RhombusC
\move(8.66025 6)
\RhombusC \RhombusC \RhombusC \RhombusC \RhombusC \RhombusC
\move(11.2583 3.5)
\RhombusC \RhombusC \RhombusC 
\move(8.66025 3)
\RhombusC \RhombusC \RhombusC \RhombusC \RhombusC 
\move(8.66025 3)
\RhombusC \RhombusC \RhombusC \RhombusC \RhombusC 
%\move(10.3923 -1)
%\RhombusC \RhombusC \RhombusC 
\move(10.3923 -1)
\RhombusB \RhombusA \RhombusB \RhombusB \RhombusB 
\move(11.2583 -.5)
\RhombusA \RhombusB \RhombusB \RhombusB \RhombusB 
%\ell
\linewd.08
\lpatt(.05 .23)
\move(6.9282 9)
\rlvec(0 -18) 
\htext(7.1 8.5){$\ell$}
\rtext td:240 (3.35 -6.0){$\sideset {} \and {} \to 
    {\left.\vbox{\vskip2.4cm}\right\}}$}
\rtext td:120 (3.44 6.15){$\sideset {} \and {} \to 
    {\left.\vbox{\vskip2.4cm}\right\}}$}
\rtext td:0 (-1.3 -.2){$\sideset {2x} \and {}\to 
    {\left\{\vbox{\vskip2.4cm}\right.}$}
\htext (2.8 6.7){$2n$}
\htext (2.75 -7.0){$2n$}
}
\centerline{\eightpoint A symmetric lozenge tiling of the hexagon $H_{n,x}$ with two holes.}
\vskip8pt
\centerline{\eightpoint Figure \FE}
}
\vskip10pt
\endinsert

\midinsert
\vskip10pt
\vbox{\noindent
\centertexdraw{
\drawdim truecm \setunitscale.7
\linewd.15
\move(0 4)
\RhombusB \RhombusB \RhombusB 
\RhombusA \RhombusB \RhombusA \RhombusB \RhombusA \RhombusA \RhombusB 
\RhombusA \RhombusB \RhombusB \RhombusA \RhombusA \RhombusA
\move(0.866025 4.5)
\RhombusB \RhombusB \RhombusB \RhombusA \RhombusB
\move(1.73205 5)
\RhombusB \RhombusB \RhombusB \RhombusA \RhombusB
\move(2.59808 5.5)
\RhombusB \RhombusB \RhombusA \RhombusB \RhombusB \RhombusA \RhombusB \RhombusB
\RhombusA \RhombusB \RhombusA \RhombusB \RhombusA 
\move(3.4641 6)
\RhombusB \RhombusB \RhombusA \RhombusB \RhombusB \RhombusA \RhombusA \RhombusB
\RhombusB \RhombusB \RhombusA 
\move(4.33013 6.5)
\RhombusB \RhombusB \RhombusA \RhombusB \RhombusA \RhombusB
\RhombusA
\move(5.19615 7)
\RhombusB \RhombusA \RhombusB \RhombusB \RhombusA
\move(6.06218 7.5)
\RhombusB \RhombusA
\move(0 1)
\RhombusC 
\move(0 0)
\RhombusC \RhombusC 
\move(0 -1)
\RhombusC \RhombusC \RhombusC \RhombusC 
\move(0 -2)
\RhombusC \RhombusC \RhombusC \RhombusC \RhombusC 
\move(0 -3)
\RhombusC \RhombusC \RhombusC \RhombusC \RhombusC 
\move(3.4641 1)
\RhombusC 
\move(3.4641 0)
\RhombusC 
\move(2.59808 -.5)
\RhombusC \RhombusC \RhombusC 
\move(4.33013 -2.5)
\RhombusC \RhombusC 
\move(5.19615 -4)
\RhombusC \RhombusC 
\move(5.19615 -5)
\RhombusC \RhombusC 
\move(5.19615 1)
\RhombusC 
%Pfade
\linewd.08
\move(6.49519 6.25)
\odSchritt \vdSchritt \odSchritt \odSchritt \odSchritt \vdSchritt
\odSchritt \odSchritt \odSchritt 
\vdSchritt \vdSchritt \vdSchritt \vdSchritt \vdSchritt
\move(6.49519 3.25)
\odSchritt \vdSchritt \odSchritt \odSchritt 
\vdSchritt \vdSchritt \vdSchritt
\odSchritt \vdSchritt \vdSchritt \vdSchritt
\move(6.49519 1.25)
 \odSchritt \vdSchritt \vdSchritt \odSchritt 
\vdSchritt \vdSchritt 
\odSchritt \vdSchritt \vdSchritt
\move(6.49519 -2.75)
\odSchritt \vdSchritt \vdSchritt \vdSchritt 
\odSchritt 
\move(6.49519 -4.75)
\vdSchritt \vdSchritt \odSchritt 
\move(2.16506 -.25)
\odSchritt \vdSchritt \vdSchritt \vdSchritt \vdSchritt 
\move(3.03109 -.75)
\vdSchritt \odSchritt \vdSchritt \vdSchritt \vdSchritt 
\ringerl(6.49519 6.25)
\ringerl(6.49519 3.25)
\ringerl(6.49519 1.25)
\ringerl(6.49519 -2.75)
\ringerl(6.49519 -4.75)
\ringerl(2.16506 -.25)
\ringerl(3.03109 -.75)
\ringerl(0.433013 -4.25)
\ringerl(1.29904 -4.75)
\ringerl(2.16506 -5.25)
\ringerl(3.03109 -5.75)
\ringerl(3.89711 -6.25)
\ringerl(4.76314 -6.75)
\ringerl(5.62917 -7.25)
\ringerl(6.49519 -7.75)
%origin
\linewd.05
\Ringerl(6.9282 0)
\htext(7.3 -.18){$O$}
%second origin
\move(7.36122 -8.25)\lcir r:.15
\htext(7.65 -8.3){$O'$}
\move(6.49519 -7.75)
\rlvec(1.29904 -.75)
\Ringerl(3.4641 0)
\htext(2.8 -.18){$C$}
}
\centerline{\eightpoint A lozenge tiling of the region
$F_{n,x}\setminus\triangleleft_2(k)$; the right boundary 
is free. The}
\centerline{\eightpoint dotted lines mark paths of lozenges. 
They determine the tiling uniquely.}
\vskip8pt
\centerline{\eightpoint Figure \FF}
}
\vskip10pt
\endinsert

The purpose of this section is to find a manageable expression for
$\M(F_{n,x}\setminus\triangleleft_2(k))$ (see Lemma~\TAE\ at the end of
this section).
In this context, we will find it more convenient to think of the tilings of
$F_{n,x}\setminus\triangleleft_2(k)$ directly  
as tilings of a half hexagon with an open boundary (cf.\ Figure~\FF)
as opposed to symmetric tilings of a hexagon with two holes
(cf.\ Figure~\FE). There is a well
known bijection between lozenge tilings of 
lattice regions and families of ``paths of lozenges'' (see Figure~{\FF}),
which in turn are equivalent to families of 
non-intersecting lattice paths (see Figure~{\FB}). Its application to
our situation is illustrated 
in Figures {\FF} and {\FB}. The origin of the system of coordinates
indicated in Figure~{\FB} corresponds to the 
point $O'$ in Figure~{\FF} (note that the bottommost path of lozenges
in Figure~{\FF} is empty for the illustrated tiling; 
the corresponding lattice path in Figure~{\FB} has no steps).

By this bijection, lozenge tilings of
$F_{n,x}\setminus\triangleleft_2(k)$ are seen to be 
equinumerous with families $(P_1,P_2,\dots,P_{2n})$ of
non-intersecting lattice paths consisting of unit horizontal and
vertical steps in the positive direction, where $P_i$ runs from $A_i=(-i,i)$
to some point from the set $I\cup\{S_1,S_2\}$, $i=1,2,\dots,2n$, with
$$\align 
I&=\{(-1,s):s=1,2,\dots,2x+2n\},\\
S_1&=(-2k-1,x+n+k),\\
S_2&=(-2k-2,x+n+k+1),
\tag\AAAAAA
\endalign$$
and the additional condition that $S_1$ and $S_2$ {\it must\/} be
ending points of some paths.

%families $(P_1,P_2,\dots,P_{2n})$ of
%non-intersecting lattice paths consisting of unit horizontal and
%vertical steps in the positive direction, where $P_i$ runs from $A_i=(-i,i)$
%to some point from the set $I\cup\{S_1,S_2\}$, $i=1,2,\dots,2n$, where
%$$\align 
%I&=\{(-1,s):s=1,2,\dots,2x+2n\},\\
%S_1&=(-2k-1,x+n+k),\\
%S_2&=(-2k-2,x+n+k+1),
%\endalign$$
%with the additional condition that $S_1$ and $S_2$ {\it must\/} be
%end points of some paths.

\midinsert
\vbox{\noindent
$$
\Gitter(1,16)(-9,0)
\Koordinatenachsen(1,16)(-9,0)
\Pfad(-2,2),122\endPfad
\Pfad(-3,3),12221\endPfad
\Pfad(-4,4),221221221\endPfad
\Pfad(-5,5),22212221121\endPfad
\Pfad(-6,6),22212\endPfad
\Pfad(-7,7),22221\endPfad
\Pfad(-8,8),22222111211121\endPfad
\DickPunkt(-1,1)
\DickPunkt(-2,2)
\DickPunkt(-3,3)
\DickPunkt(-4,4)
\DickPunkt(-5,5)
\DickPunkt(-6,6)
\DickPunkt(-7,7)
\DickPunkt(-8,8)
\DickPunkt(-1,4)
\DickPunkt(-1,6)
\DickPunkt(-1,10)
\DickPunkt(-1,12)
\DickPunkt(-1,15)
\DickPunkt(-5,10)
\DickPunkt(-6,11)
\Kreis(0,0)
\Label\u{A_1}(-1,1)
\Label\u{A_2}(-2,2)
\Label\u{A_3}(-3,3)
\Label\u{A_4}(-4,4)
\Label\u{A_5}(-5,5)
\Label\u{A_6}(-6,6)
\Label\u{A_7}(-7,7)
\Label\u{A_8}(-8,8)
\Label\ro{S_1}(-5,10)
\Label\ro{S_2}(-6,11)
\Label\ro{P_1}(-1,1)
\Label\r{P_2}(-1,3)
\Label\r{P_3}(-2,5)
\Label\r{P_4}(-4,7)
\Label\ro{P_5}(-5,7)
\Label\o{P_6}(-6,9)
\Label\o{P_7}(-7,11)
\Label\l{P_8}(-8,12)
\Label\ro{O'}(0,0)
\Label\r{x}(1,0)
\Label\ru{y}(0,16)
\hskip-4cm
$$
\centerline{\eightpoint The paths of lozenges of Figure {\FF} drawn as
non-intersecting lattice paths on $\Z^2$.}
\vskip8pt
\centerline{\eightpoint Figure \FB}
}
\vskip10pt
\endinsert

At this point, we need a slight extension of Theorem~3.2 in
\cite{\StemAE} (which is, in fact, derivable from the minor summation
formula of Ishikawa and Wakayama \cite{\IsWaAA, Theorem~2}). 
The reader should recall that the {\it Pfaffian} of a skew-symmetric
$2n\times 2n$ matrix $A$ can be defined by (see e.g.\ 
\cite{\StemAE, p.~102})
$$\Pf A:=
\sum _{\pi\in\Cal M[1,\dots, 2n]} ^{}\sgn\pi 
\underset i,j\text{ matched in }\pi\to{\prod _{i<j} ^{}}
A_{i,j},
\tag\PPF$$
where $M[1,2,\dots, 2n]$ denotes the set of all perfect matchings
($1$-factors) of (the complete graph on) $\{1,2,\dots,2n\}$, and where
$\sgn\pi=(-1)^{\text{cr}(\pi)}$, with $\text{cr}(\pi)$ denoting the
number of ``crossings" of $\pi$, that is, the number of quadruples
$i<j<k<l$ such that, under $\pi$, $i$ is paired with $k$, and $j$ is
paired with $l$. It is a well-known fact
(see e.g\. \cite{\StemAE, Prop.~2.2}) that
$$
(\Pf A)^2=\det A.
\tag\PF$$

\proclaim{Theorem \TAD}
Let $\{A_1,A_2,\dots,A_{p},S_1,S_2,\dots,S_q\}$ and
$I=\{I_1,I_2,\dots\}$ be finite sets of lattice points in the integer lattice
$\Z^2$. Then 
$$
\Pf\pmatrix \hphantom{-}Q&H\\-H^t&0\endpmatrix=
(-1)^{\binom q2}\sum _{\pi\in S_{p}} ^{}(\sgn\pi)\cdot
\Cal P^{\text{nonint}}(\bold A_\pi\to \bold S\cup I),
\tag\AAAAA
$$
where $\bold A_\pi=(A_{\pi(1)},A_{\pi(2)},\dots,A_{\pi(p)})$,
and $\Cal P^{\text{nonint}}(\bold A_\pi\to \bold S\cup I)$ is
the number of families $(P_1,P_2,\dots,P_{p})$ of
non-intersecting lattice paths consisting of unit horizontal and
vertical steps in the positive direction, with $P_k$ running from
$A_{\pi(k)}$ to $S_k$, for $k=1,2,\dots,q$, and to $I_{j_k}$, for
$k=q+1,q+2,\dots,p$, the indices being required to satisfy 
$j_{q+1}<j_{q+2}<\dots<j_p$. The matrix $Q=(Q_{i,j})_{1\le i,j\le p}$
is defined by 
$$Q_{i,j}=
\sum _{1\le s<t} ^{}
\big(\Cal P(A_i\to I_s)\cdot \Cal P(A_j\to I_t)-
\Cal P(A_j\to I_s)\cdot \Cal P(A_i\to I_t)\big),
\tag\AAAAAAA
$$
where $\Cal P(A\to E)$ denotes the number of lattice paths from $A$ to
$E$, and the matrix $H=(H_{i,j})_{1\le i\le p,\ 1\le j\le q}$ by 
$$H_{i,j}=\Cal P(A_i\to S_{j}).$$
\endproclaim

In the special case when the starting and ending points satisfy a
certain compatibility condition (called $D$-compatibility in
\cite{\StemAE}), the only permutation $\pi$ which contributes to the
right-hand side of (\AAAAA) is the identity permutation, and (\AAAAA)
reduces to \cite{\StemAE, Theorem~3.2}. In our context, the
compatibility condition is not satisfied. However, the same arguments that
prove \cite{\StemAE, Theorem~3.2} can be used to obtain (\AAAAA).
(Alternatively, one could use the minor summation formula
of Ishikawa and Wakayama \cite{\IsWaAA, Theorem~2}. In it, choose 
$m=p$, $r=q$, and the skew-symmetric matrix $B$ to be $B_{i,j}=1$ for $i<j$
--- which makes all principal Pfaffian minors of $B$ equal $1$ --- to expand
the Pfaffian on the left-hand side of (\AAAAA) into a sum of minors of
a certain matrix. Each minor can then be seen to count certain families of
nonintersecting lattice paths by the general form of the
Lindstr\"om--Gessel--Viennot theorem \cite{\LindAA, Lemma~1},
\cite{\GeViAB, Theorem~1}, and, altogether, these are the families that
are described in the statement of Theorem~\TAD.)

We now apply Theorem~\TAD\ to our situation, that is, $p=2n$, $q=2$,
$A_i=(-i,i)$, for $i=1,2,\dots,2n$, and $S_1$, $S_2$, and $I$ are
given by (\AAAAAA). It is not difficult to convince oneself that, for
this choice of starting and ending points, all families of
nonintersecting lattice paths counted on the right-hand side of
(\AAAAA) give rise to even permutations $\pi$. Hence, the right-hand
side of (\AAAAA) counts indeed the families of nonintersecting lattice
paths that we need to count. By Theorem~\TAD, their number 
is equal to the negative value of the Pfaffian of
$$M_n(x):=\pmatrix \hphantom{-}Q&H\\-H^t&0\endpmatrix,\tag\AAAA$$
where $Q$ is a $(2n)\times(2n)$ skew-symmetric matrix with
$(i,j)$-entry $Q_{i,j}$ given by (\AAAAAAA), 
and where $H$ is a $(2n)\times2$ matrix,
in which the $(i,j)$-entry $H_{i,j}$ is equal to the number of paths from
$A_i$ to $S_j$, $i=1,2,\dots,2n$, $j=1,2$. 
(It is the negative value
of the Pfaffian because of the sign $(-1)^{\binom q2}$ on the
right-hand side of (\AAAAA), as we have $q=2$.)

In particular, using the fact that the number of lattice paths on the
integer lattice $\Z^2$ between
two given lattice points is given by a binomial coefficient, we have
$$\align 
H_{i,1}&=\binom {x+n-k-1}{i-2k-1},\tag\AAAa\\
H_{i,2}&=\binom {x+n-k-1}{i-2k-2}.\tag\AAAb
\endalign$$
On the other hand,
substituting $A_i=(-i,i)$ and $I_s=(-1,s)$ in (\AAAAAAA),
we have
$$\align 
Q_{i,j}&=\sum _{1\le s<t\le 2x+2n} ^{}(\vert\Cal P(A_i\to I_s)\vert\cdot\vert
\Cal P(A_j\to I_t)\vert
-\vert\Cal P(A_j\to I_s)\vert\cdot\vert \Cal P(A_i\to I_t)\vert)\\
&=\sum _{1\le s<t\le 2x+2n} ^{}\bigg(\binom {s-1}{i-1}\binom
{t-1}{j-1}-\binom {s-1}{j-1}\binom
{t-1}{i-1}\bigg)\\
&=\sum _{1\le s\le t\le 2x+2n} ^{}\binom {s-1}{i-1}\binom
{t-1}{j-1}-\sum _{1\le s\le t\le 2x+2n} ^{}\binom {s-1}{j-1}\binom
{t-1}{i-1}\\
&=\sum _{1\le t\le 2x+2n} ^{}\binom {t}{i}\binom
{t-1}{j-1}-\sum _{1\le t\le 2x+2n} ^{}\binom {t}{j}\binom
{t-1}{i-1}\tag\AAa\\
&=\sum _{t=1} ^{2x+2n}\frac {j-i} {t}\binom {t}{i}\binom
{t}{j},
\tag\AA
\endalign$$
where we used the well-known identity
$$
\sum _{t=0} ^{X}\binom t{i-1}=\binom {X+1}i
\tag\AAb$$
to obtain (\AAa).
We may obtain an alternative expression for $Q_{i,j}$ by replacing 
$\frac {1} {t}\binom ti=\frac {1} {i}\binom {t-1}{i-1}$ in the last
expression by $\frac {1} {i}\sum _{l=0} ^{i-1}\binom {t-j}l\binom
{j-1}{i-l-1}$, this equality being true because of the Chu--Vandermonde
summation (cf\. e.g\. \cite{\GrKPAA, Sec.~5.1, (5.27)}).
Thus, we arrive at
$$\align Q_{i,j}
&=\sum _{l=0} ^{i-1}\sum _{t=1} ^{2x+2n}\frac {j-i} {i}\binom {t-j}l\binom
{j-1}{i-l-1}\binom {t}{j}\\
&=\sum _{l=0} ^{i-1}\sum _{t=1} ^{2x+2n}\frac {j-i} {i}\binom
{j-1}{i-l-1}\binom {l+j}l
\binom {t}{l+j}\\
&=\sum _{l=0} ^{i-1}\frac {j-i} {i}\binom
{j-1}{i-l-1}\binom {l+j}l
\binom {2x+2n+1}{l+j+1},
\tag\AB
\endalign$$
the last line again being due to (\AAb).

To summarize, we have obtained the following result.

\proclaim{Lemma \TAE}
For all positive integers $n,x$ and nonnegative integers $k$, we have
$$
\M(F_{n,x}\setminus\triangleleft_2(k))=
-\Pf M_n(x),
\tag\ABb$$
where $M_n(x)$ is given by {\rm(\AAAA)}, with
$Q_{i,j}$ defined in {\rm(\AA)} or {\rm(\AB)}, and $H_{i,j}$ defined in {\rm(\AAAa)}
and {\rm(\AAAb)}.
\endproclaim

\head 5. Proof of Theorem {\TB}\endhead

In the sequel, 
we shall interpret sums by
$$\sum _{k=m} ^{n-1}\text {\rm Expr}(k)=\cases \hphantom{-}
\sum _{k=m} ^{n-1} \text {\rm Expr}(k)&n>m\\
\hphantom{-}0&n=m\\
-\sum _{k=n} ^{m-1}\text {\rm Expr}(k)&n<m.\endcases
\tag\ABB$$
In particular, using this convention, the expression for
$Q_{i,j}$ given in (\AA) makes sense for negative integers
$x$ also (in which case the upper bound in the sum can be negative) 
and is actually equal to the expression in (\AB). 
It is the latter fact that we shall frequently make use of.

\medskip
Our proof of Theorem~\TB\ involves a sequence of five steps.
By Lemma~\TAE, we know that the number that we want to compute 
is the negative of
a Pfaffian. We shall frequently use the fact (\PF)
that the square of the Pfaffian of a
skew-symmetric matrix is equal to its determinant. 

By its definition, $\Pf M_n(x)$ is a polynomial in $x$.
In Step~1 we prove that
$$\det M_n(x)=\det M_n(-2n-x ).$$
With $d$ denoting the degree of $\Pf M_n(x)$ as a polynomial in $x$,
this implies that
$$\Pf M_n(x)=(-1)^d\Pf M_n(-2n-x ).
\tag\SS$$
Subsequently, in Step~2 we show that
$$\underset s\ne n-k\to{\prod _{s=1} ^{n}}(x+s)_{2n-2s+1}^2$$
divides $\det M_n(x)$ as a polynomial in $x$, while in Step~3 we show
that 
$${\prod _{s=1} ^{n-1}}(x+s+\tfrac {1} {2})_{2n-2s}^2$$
divides $\det M_n(x)$. Both combined, this proves that
$$\underset s\ne n-k\to{\prod _{s=1} ^{n}}(x+s)_{2n-2s+1}
{\prod _{s=1} ^{n-1}}(x+s+\tfrac {1} {2})_{2n-2s},$$
which is a polynomial of degree 
$$\sum _{s=1} ^{n}(4n-4s+1)-(2k+1)=n(2n-1)-(2k+1),$$
divides $\Pf M_n(x)$ as a polynomial in $x$. The computation in
Step~4 then shows that the degree of $\Pf M_n(x)$, as a polynomial in
$x$, is at most $2n^2+n-4k-3$. Altogether, this implies that
$$-\Pf M_n(x)=P_n(x)
\underset s\ne n-k\to{\prod _{s=1} ^{n}}(x+s)_{2n-2s+1}
{\prod _{s=1} ^{n}}(x+s+\tfrac {1} {2})_{2n-2s},
\tag\BD$$
where $P_n(x)$ is a polynomial in $x$ of degree at most 
$$2n^2+n-4k-3-n(2n-1)+(2k+1)=2n-2k-2.$$ 

In Step~5, we determine the value of $P_n(x)$ at
$x=0,-1,\dots,-n+k+1$ (see (\BI)). The corresponding calculations make use of an
auxiliary lemma due to Mehta and Wang \cite{\MeWaAA}, see Theorem~\TAF\
and Corollary~\TCb\ in Section~6.
By (\SS), this
gives us at the same time the value of $P_n(x)$ at
$x=-2n,-2n+1,\dots,-n-k-1$. In total, these are $2n-2k$ explicit 
evaluations of $P_n(x)$ at special values of $x$. Given the fact that
the degree of $P_n(x)$ is at most $2n-2k-2$, they determine $P_n(x)$
uniquely, and an explicit expression for $P_n(x)$ can be written down
using Lagrange interpolation. If this is substituted into (\BD), then
the evaluation of $-\Pf M_n(x)$ is complete. After some manipulations,
one arrives at the expression in (\ABa).

\smallskip
{\smc Step 1. $\det M_n(x)=\det M_n(-2n-x )$}.
We prove this claim by transforming, up to sign,
$M_n(x)$ into $M_n(-2n-x)$ by a sequence of elementary row and column
operations (which, of course, leave the value of the determinant invariant).
To be precise, for $i=2n,2n-1,\dots,2$ (in this order), we add
$$\sum _{a=1} ^{i-1}\binom {i-1}{a-1}\cdot(\text {row $a$})$$
to row $i$, and then for $j=2n,2n-1,\dots,2$, we add
$$\sum _{b=1} ^{j-1}\binom {j-1}{b-1}\cdot(\text {column $b$})$$
to column $j$. Let $M^{(1)}_n(x)$ denote the matrix which arises
after these row and column operations. 
According to (\AAa), the $(i,j)$-entry in $M^{(1)}_n(x)$ is
$$
\sum_{a = 1}^{i}
\sum_{b = 1}^{j}\binom {i-1}{a-1}\binom {j-1}{b-1}
\sum _{t=1} ^{2x+2n}
\left(\binom {t}{a}\binom
{t-1}{b-1}-\binom {t-1}{a-1}\binom {t}{b}\right)
\tag\AC$$
for $1\le i,j\le 2n$. By (\AAAa) and (\AAAb), 
for $1\le i\le 2n$ and $j=2n+\ep$, $\ep=1,2$, the $(i,j)$-entry of 
$M^{(1)}_n(x)$ is
$$
\sum_{a = 1}^{i}\binom {i-1}{a-1}
\binom {x+n-k-1}{a-2k-\ep},
\tag\AD$$
and, for $1\le j\le 2n$ and $i=2n+\ep$, $\ep=1,2$, it is
$$
-\sum_{b = 1}^{j}\binom {j-1}{b-1}
\binom {x+n-k-1}{b-2k-\ep}.
\tag\AE$$
By Chu--Vandermonde summation, we have
$$
\sum_{a = 1}^{i}\binom {i-1}{a-1}\binom {t+\gamma}{a+\eta}=
\sum_{a = 1}^{i}\binom {i-1}{i-a}\binom {t+\gamma}{a+\eta}=
\binom {t+i+\gamma-1}{i+\eta},
$$
whence the expression (\AC) simplifies to
$$\align 
\sum _{t=1} ^{2x+2n}&
\left(\binom {t+i-1}{i}\binom
{t+j-2}{j-1}-\binom {t+i-2}{i-1}\binom {t+j-1}{j}\right)\\
&=\sum _{t=-2x-2n} ^{-1}
\left(\binom {-t+i-1}{i}\binom
{-t+j-2}{j-1}-\binom {-t+i-2}{i-1}\binom {-t+j-1}{j}\right)\\
&=(-1)^{i+j-1}\sum _{t=-2x-2n} ^{-1}
\left(\binom {t}{i}\binom
{t}{j-1}-\binom {t}{i-1}\binom {t}{j}\right)\\
&=(-1)^{i+j-1}\sum _{t=-2x-2n} ^{-1}
\left(\binom {t+1}{i}\binom
{t}{j-1}-\binom {t}{i-1}\binom {t+1}{j}\right)
\tag\AFa
\\
&=(-1)^{i+j}\sum _{t=0} ^{-2x-2n-1}
\left(\binom {t+1}{i}\binom
{t}{j-1}-\binom {t}{i-1}\binom {t+1}{j}\right)
\tag\AFb
\\
&=(-1)^{i+j}\sum _{t=1} ^{-2x-2n}
\left(\binom {t}{i}\binom
{t-1}{j-1}-\binom {t-1}{i-1}\binom {t}{j}\right).
\tag\AF
\endalign$$
Here, we used the identity 
$$\align 
\binom {t}{i}\binom
{t}{j-1}-\binom {t}{i-1}\binom {t}{j}&=
\(\binom {t}{i}+\binom {t}{i-1}\)\binom
{t}{j-1}-\binom {t}{i-1}
\(\binom {t}{j}+\binom {t}{j-1}\)\\
&=
\binom {t+1}{i}\binom
{t}{j-1}-\binom {t}{i-1}\binom {t+1}{j}
\endalign$$
to obtain (\AFa),
and our convention (\ABB) for sums to obtain (\AFb).
Comparison with (\AAa) shows that this last expression is, up to the
sign $(-1)^{i+j}$, exactly
$Q_{i,j}$ with $x$ replaced by $-2n-x$. In a similar vein, the
expression (\AD) simplifies to
$$
\binom {x+n+i-k-2}{i-2k-\ep}=(-1)^{i+\ep}\binom
{-x-n-k+1-\ep}{i-2k-\ep},\tag\AG$$
while expression (\AE) simplifies to the same expression with $i$
replaced by $j$. Upon setting $\ep=2$, this shows that 
the $(i,2n+2)$-entry in $M^{(1)}_n(x)$ is, up
to the sign $(-1)^{i}$, identical with the $(i,2n+2)$-entry in
$M_n(-2n-x)$, with an analogous statement being true for the
$(2n+2,j)$-entry of $M^{(1)}_n(x)$ and the $(2n+2,j)$-entry of
$M_n(-2n-x)$. 

We do one last row and one last column operation: in $M^{(1)}_n(x)$,
we add the last row to the next-to-last row, and we add 
the last column to the next-to-last column. Let $M^{(2)}_n(x)$ denote
the resulting matrix. By (\AG),
for $i=1,2,\dots,2n$, the $(i,2n+1)$-entry of
$M^{(2)}_n(x)$ is equal to
$$
(-1)^{i+1}\binom {-x-n-k}{i-2k-1}+
(-1)^{i}\binom {-x-n-k-1}{i-2k-2}=
(-1)^{i+1}\binom {-x-n-k-1}{i-2k-1},
\tag\AH$$
which is, up to the sign $(-1)^{i+1}$, exactly the $(i,2n+1)$-entry
in $M_n(-2n-x)$. An analogous statement is true for the
$(2n+1,j)$-entries of $M^{(2)}_n(x)$ and $M_n(-2n-x)$.

In summary, 
as the two by two block in the lower right corner of $M_n(x)$ consists of 
zeros,
the computations (\AF)--(\AH) show that the $(i,j)$-entry
of $M^{(2)}_n(x)$ is $(-1)^{i+j}$ times the $(i,j)$-entry of
$M_n(-2n-x)$. This implies
$$\det M_n(x)=\det M^{(2)}_n(x)=\det M_n(-2n-x),$$
as claimed.

\smallskip
{\smc Step 2}.
{\it $\dsize
\underset s\ne n-k\to{\prod _{s=1} ^{n}}(x+s)_{2n-2s+1}^2$ divides $\det
M_n(x)$}. 
We begin by observing that the product in the claim can be
also rewritten as
$$\underset s\ne n-k\to{\prod _{s=1} ^{n}}(x+s)_{2n-2s+1}^2=
\prod _{s=1} ^{n}(x+s)^{2s+2\chi(s<n-k)-2}
\prod _{s=n+1} ^{2n-1}(x+s)^{4n-2s+2\chi(s>n+k)-2},
\tag\BAa
$$
where $\chi(\Cal A)=1$ if $\Cal A$ is
true and $\chi(\Cal A)=0$ otherwise. In view of Step~1, it suffices
to establish that 
$$\prod _{s=1} ^{n}(x+s)^{2s+2\chi(s<n-k)-2}
$$
divides $\det M_n(x)$. 

Now let $s$, $a$ and $b$ be integers with $1\le s\le n$ and $1\le a\le b\le
2n$. We claim that
$$\sum _{i=a} ^{b}\binom {b-a}{i-a}(\text {row $i$ of $M_n(-s)$})=0,
\tag\BA$$
as long as 
\roster
\item"(A)" $a-b\le 2n-2s<a$, and
\item"(B)" either $b\le 2k$, or $b\ge 2k+3$ and $a-b+k+1\le n-s<a-k-1$.
\endroster

Indeed, if we specialize (\BA) to the $j$-th column,
where $j\le 2n$, we obtain, using the expression (\AA) for $Q_{i,j}$,
$$\align 
\sum _{i=a} ^{b}\binom {b-a}{i-a}Q_{i,j}&=
\sum _{i=a} ^{b}\binom {b-a}{i-a}
\sum _{t=1} ^{2n-2s}\frac {j-i} {t}\binom {t}{i}\binom
{t}{j}\\
&=
\sum _{t=1} ^{2n-2s}\sum _{i=a} ^{b}\binom {b-a}{b-i}
\left(\binom {t}{i}\binom {t-1}{j-1}-\binom {t-1}{i-1}\binom {t}{j}\right)\\
&=
\sum _{t=1} ^{2n-2s}\left(\binom {b-a+t}{b}
\binom {t-1}{j-1}-\binom {b-a+t-1}{b-1}\binom {t}{j}\right).
\tag\BB
%\\&=
%-\sum _{t=2n-2s+1} ^{0}\left(\binom {b-a+t}{b}
%\binom {t-1}{j-1}-\binom {b-a+t-1}{b-1}\binom {t}{j}\right).
%\tag\BC
\endalign$$
Here we used Chu--Vandermonde summation (cf.\ 
\cite{\GrKPAA, Sec.~5.1, (5.27)}) in the last line.
Since, for $a-b\le 2n-2s<a$ (which is condition (A)), 
the binomial coefficients containing the 
parameter $b$ in expression (\BB) 
are identically zero throughout the summation range,
it is clear that the corresponding sum vanishes.

On the other hand, if we specialize (\BA) to the $(2n+1)$-st column,
we obtain, again using Chu--Vandermonde summation,
$$\align 
\sum _{i=a} ^{b}\binom {b-a}{i-a}H_{i,1}\bigg\vert_{x=-s}&=
\sum _{i=a} ^{b}\binom {b-a}{b-i}\binom {n-s-k-1}{i-2k-1}\\
&=
\binom {b-a+n-s-k-1} {b-2k-1},
\endalign$$
which vanishes for $b\le 2k$, and for $b\ge2k+2$ and 
$0\le b-a+n-s-k-1< b-2k-1$, the last inequality being equivalent to
$a-b+k+1\le n-s< a-k$, and
if we specialize (\BA) to the $(2n+2)$-nd column,
we obtain
$$\align 
\sum _{i=a} ^{b}\binom {b-a}{i-a}H_{i,2}\bigg\vert_{x=-s}&=
\sum _{i=a} ^{b}\binom {b-a}{b-i}\binom {n-s-k-1}{i-2k-2}\\
&=
\binom {b-a+n-s-k-1} {b-2k-2},
\endalign$$
which vanishes for  $b\le 2k-1$, and for $b\ge2k+3$ and
$0\le b-a+n-s-k-1< b-2k-2$, the last inequality being equivalent to
$a-b+k+1\le n-s< a-k-1$.
This establishes our claim.

In order to prove that $(x+s)^{2s}$ divides $\det M_n(x)$ for $1\le
s<n-k$, we use (\BA) with $a=2n-2s+1$ and $2n-2s+1\le b\le 2n$.
It is not difficult to see that for these choices of $s$, $a$ 
and $b$ the
conditions (A) and (B) are satisfied, so that we obtain $2s$ linear
combinations of the rows that are linearly independent 
(as, for our choices of $a$ and $b$, the coefficients in (\BA) form a triangular 
array) and vanish when $x=-s$. This implies
divisibility by $(x+s)^{2s}$ (cf\. e.g\. \cite{\KratBI, Lemma in Sec.~2}).

To prove that $(x+s)^{2s-2}$ divides $\det M_n(x)$ for $n-k\le
s\le n$, we use (\BA) with $a=2n-2s+1$ and $2n-2s+1\le b\le 2k$ on
the one hand, and with $a=n+k-s+2$ and $2k+3\le b\le 2n$ on the other
hand.
Again, it is not difficult to see that for both types of choices of 
$s$, $a$ and $b$ 
the conditions (A) and (B) are satisfied, so that we obtain
$(2k+2s-2n)+(2n-2k-2)=2s-2$ linear
combinations of the rows that are linearly independent and vanish
at $x=-s$. In the same way as before, 
this implies divisibility by $(x+s)^{2s-2}$.

%In order to prove that $(x+s)^{4n-2s-2}$ divides $\det M_n(x)$ for $n\le
%s\le n+k$, we use (\BA) with $a=1$ and $2s-2n+1\le b\le 2k$ on
%the one hand, and with $a=n+k-s+2$ and $2k+3\le b\le 2n$ on the other
%hand.
%Again, it is not difficult to see that for both types of choices of $s,a,b$ 
%the conditions (A) and (B) are satisfied, so that we obtain
%$(2k+2n-2s)+(2n-2k-2)=4n-2s-2$ linear
%combinations of the rows that are linearly independent, as required.
%
%Finally, in order to prove that $(x+s)^{4n-2s}$ divides $\det M_n(x)$ for 
%$n+k<s\le 2n$, we use (\BA) with $a=1$ and $2s-2n+1\le b\le 2n$.
%Also here, it is not difficult to see that for these choices of $s,a,b$ the
%conditions (A) and (B) are satisfied, so that we obtain $4n-2s$ linear
%combinations of the rows that are linearly independent, as required.

\smallskip
{\smc Step 3}.
{\it $\dsize
{\prod _{s=1} ^{n-1}}(x+s+\tfrac {1} {2})_{2n-2s}^2$ divides $\det
M_n(x)$}. 
We begin by observing that the product in the claim can be
also rewritten as
$${\prod _{s=1} ^{n-1}}(x+s+\tfrac {1} {2})_{2n-2s}^2=
\prod _{s=1} ^{n-1}(x+s+\tfrac {1} {2})^{2s}
\prod _{s=n} ^{2n-2}(x+s+\tfrac {1} {2})^{4n-2s-2}.
$$
In view of Step~1, it suffices
to establish that 
$$\prod _{s=1} ^{n-1}(x+s+\tfrac {1} {2})^{2s}
$$
divides $\det M_n(x)$.

In order to prove the claim for $s< n-k$, we shall show that
$$\multline 
(n-s-k)\cdot\left(\text {row
$(2n-2s+1)$ of $M_n(-s-\tfrac {1} {2})$}\right)\\+
\frac {1} {2}(n-s-k)\cdot\(\text {row
$(2n-2s)$ of $M_n(-s-\tfrac {1} {2})$}\)\\
+
\sum _{i=1} ^{2n-2s-1}\frac {(-1)^{i}} {2^{2n-2s-i+2}}\(\text {row
$i$ of $M_n(-s-\tfrac {1} {2})$}\)=0,
\endmultline\tag\BCa$$
and that
$$\multline 
\(\text {row $i$ of $M_n(-s-\tfrac {1} {2})$}\)+
\frac {2i+2s-2n-2k-3} {i-2k-1}
\(\text {row $(i-1)$ of $M_n(-s-\tfrac {1} {2})$}\)\\
+
\frac {(2i+2s-2n-2k-3)(2i+2s-2n-2k-5)} {4(i-2k-1)(i-2k-2)}
\(\text {row $(i-2)$ of $M_n(-s-\tfrac {1} {2})$}\)
=0
\endmultline\tag\BCb$$
for $i=2n-2s+2,2n-2s+3,\dots,2n$. As these are linearly independent row
combinations, the claim will follow.

In order to prove the claim for $s\ge n-k$, we shall show that
$$\sum _{i=1} ^{2n-2s-1}\frac {(-1)^{i}} {2^{2n-2s-i-1}}\(\text {row
$i$ of $M_n(-s-\tfrac {1} {2})$}\)=0,\tag\BCc$$
that 
$$\(\text {row $i$ of $M_n(-s-\tfrac {1} {2})$}\)=0\tag\BCd$$
for $i=2n-2s,2n-2s+1,\dots,2k$, and that 
$$\multline 
\(\text {row $i$ of $M_n(-s-\tfrac {1} {2})$}\)+
\frac {2i+2s-2n-2k-3} {i-2k-1}
\(\text {row $(i-1)$ of $M_n(-s-\tfrac {1} {2})$}\)\\
+
\frac {(2i+2s-2n-2k-3)(2i+2s-2n-2k-5)} {4(i-2k-1)(i-2k-2)}
\(\text {row $(i-2)$ of $M_n(-s-\tfrac {1} {2})$}\)
=0
\endmultline\tag\BCe$$
for $i=2k+3,2k+4,\dots,2n$.
Again, as these are linearly independent row
combinations, the claim will follow.

Let first $s\ge n-k$.
We start with the proof of (\BCc). Specializing (\BCc) to the $j$-th
column, $j=1,2,\dots,2n$, by (\AAa) we see that we must prove the
identity
$$\sum _{i=1} ^{2n-2s-1}\frac {(-1)^{i}} {2^{2n-2s-i-1}}
\Bigg(\sum _{t=1} ^{2n-2s-1} \binom {t}{i}\binom
{t-1}{j-1}-\sum _{t=1} ^{2n-2s-1} \binom {t-1}{i-1}
\binom {t}{j}\Bigg)=0.
\tag\BCf$$
In order to see that this is indeed true, we first extend the sum
over $i$ to the range $i=0,1,\dots,2n-2s-1$, thereby obtaining
$$\align 
\sum _{i=0} ^{2n-2s-1}&\frac {(-1)^{i}} {2^{2n-2s-i-1}}
\Bigg(\sum _{t=1} ^{2n-2s-1} \binom {t}{i}\binom
{t-1}{j-1}-\sum _{t=1} ^{2n-2s-1} \binom {t-1}{i-1}
\binom {t}{j}\Bigg)\\
&\kern7cm-
\frac {1} {2^{2n-2s-1}}\sum _{t=1} ^{2n-2s-1} \binom {t-1}{j-1}\\
&=\sum _{i=0} ^{2n-2s-1}\frac {(-1)^{i}} {2^{2n-2s-i-1}}
\Bigg(\sum _{t=1} ^{2n-2s-1} \binom {t}{i}\binom
{t-1}{j-1}-\sum _{t=1} ^{2n-2s-1} \binom {t-1}{i-1}
\binom {t}{j}\Bigg)\\
&\kern7cm-\frac {1} {2^{2n-2s-1}} \binom {2n-2s-1}{j}
\endalign$$
for the left-hand side of (\BCf). Next we interchange the
sum over $i$ with the sums over $t$, and subsequently we evaluate the
(now inner) sums over $i$ by means of the binomial theorem. In this
manner, the left-hand side of (\BCf) becomes
$$\align 
\frac {1} {2^{2n-2s-1}}&
\sum _{t=1} ^{2n-2s-1} (1-2)^t\binom
{t-1}{j-1}+\frac {1} {2^{2n-2s-2}}
\sum _{t=1} ^{2n-2s-1} (1-2)^{t-1}\binom {t}{j}\\
&\kern7cm
-\frac {1} {2^{2n-2s-1}} \binom {2n-2s-1}{j}\\
&=\frac {1} {2^{2n-2s-1}}
\sum _{t=1} ^{2n-2s-1} (-1)^t\frac {(t-1)!} {j!\,(t-j)!}(j-2t)
-\frac {1} {2^{2n-2s-1}} \binom {2n-2s-1}{j}\\
&=-\frac {1} {2^{2n-2s-1}}\Bigg(
\sum _{t=1} ^{2n-2s-1} (-1)^t\Bigg(\binom tj+\binom {t-1}j\Bigg)
+\binom {2n-2s-1}{j}\Bigg)=0,
\endalign$$
as desired.

On the other hand, specializing (\BCc) to the $j$-th
column, $j=2n+1,2n+2$, by (\AAAa) and (\AAAb) we obtain
$$\sum _{i=1} ^{2n-2s-1}\frac {(-1)^{i}} {2^{2n-2s-i-1}}
\binom {n-k-s-\frac {3} {2}}{i-2k-\ep},
$$
where $\ep=1,2$, which is indeed zero since
the binomial coefficient always vanishes because of $i\le
2n-2s-1<2k+\ep$, the last inequality being due to our assumption $s\ge n-k$.

That (\BCd) holds can be easily checked by inspection.

For the proof of (\BCe), we observe that we have 
$Q_{i,j}\big\vert_{x=-s-\frac {1} {2}}=0$
for all $i\ge 2n-2s$, because in this case the appearance of the
binomial coefficient $\binom ti$ in the sum in formula (\AA)
implies that all summands of this sum vanish. In its turn, this
entails that the left-hand side of (\BCe) specialized to the $j$-th
column, where $1\le j\le 2n$, is trivially zero since
$$i>i-1>i-2\ge 2k+1\ge 2n-2s+1>2n-2s,$$
by our assumptions. To see that the left-hand side of (\BCe) is as well zero
when it is specialized to the $(2n+1)$-st or $(2n+2)$-nd column
amounts to a routine verification using the expressions (\AAAa) and
(\AAAb) for the corresponding matrix entries.

We now assume that $s<n-k$ and
turn our attention to (\BCa). The reader should notice that the relations 
(\BCa) and (\BCc) are relatively similar, the essential difference being
the two extra terms in (\BCa) corresponding to the $(2n-2s)$-th and the
$(2n-2s+1)$-st row, respectively. If $1\le j\le 2n$,
the proof of relation (\BCa) specialized
to column $j$ is therefore identical with the proof of relation (\BCc)
specialized to column $j$, because the entries in the first $2n$ columns
of the $(2n-2s)$-th and the $(2n-2s+1)$-st row evaluated at
$x=-s-\frac {1} {2}$ are all zero. (The reader should recall formula
(\AA).) To show the relation (\BCa) specialized to the $(2n+1)$-st
respectively to the $(2n+2)$-nd column requires however more work.
We have to prove
$$\multline 
(n-s-k)\left(\binom {n-s-k-\frac {3} {2}}{2n-2s-2k-\ep+1}+
\frac {1} {2}\binom {n-s-k-\frac {3} {2}}{2n-2s-2k-\ep}\right)\\
+
\sum _{i=1} ^{2n-2s-1}\frac {(-1)^{i}} {2^{2n-2s-i+2}}
\binom {n-s-k-\frac {3} {2}}{i-2k-\ep}=0,
\endmultline$$
where $\ep=1,2$, respectively, after simplification,
$$\multline 
\frac {(n-s-k)(\ep-2)}2\cdot\frac {(-n+s+k+\ep-\frac {1} {2})_{2n-2s-2k-\ep}} 
{(2n-2s-2k-\ep+1)!}\\
+
\sum _{i=1} ^{2n-2s-1}\frac {(-1)^{i}} {2^{2n-2s-i+2}}
\binom {n-s-k-\frac {3} {2}}{i-2k-\ep}=0.
\endmultline\tag\BCg$$
We reverse the order of summation in the sum over $i$ (that is, we
replace $i$ by $2n-2s-i-1$), and subsequently we write the (new)
sum over $i$ in standard hypergeometric notation
$${}_p F_q\!\left[\matrix a_1,\dots,a_p\\ b_1,\dots,b_q\endmatrix; 
z\right]=\sum _{m=0} ^{\infty}\frac {\po{a_1}{m}\cdots\po{a_p}{m}}
{m!\,\po{b_1}{m}\cdots\po{b_q}{m}} z^m\ .
\tag\HYP$$
Thereby we obtain 
$$\multline 
\frac {(n-s-k)(\ep-2)}2\cdot\frac {(-n+s+k+\ep-\frac {1} {2})_{2n-2s-2k-\ep}} 
{(2n-2s-2k-\ep+1)!}\\
-{\frac {   
       ({ \textstyle - n  + k + s+ \ep +{\frac 1 2}}) _{
         2 n - 2 s-2k-\ep-1}   }
   {8\,({2 n - 2 s-2k-\ep-1 })! }}
{} _{2} F _{1} \!\left [ \matrix {1, - 2 n + 2 k +
        2 s+ \ep +1}\\ {  -n + k + s+\ep+{\frac 1 2}}\endmatrix ;
        {\displaystyle {\frac 1 2}}\right ]
\endmultline\tag\BCh$$
for the left-hand side of (\BCg).

If $\ep=2$, then the $_2F_1$-series in (\BCh) can be evaluated using
Gau{\ss}' second $_2F_1$-summation
(cf\. \cite{\SlatAC, (1.7.1.9); Appendix (III.6)})
$${}_2F_1\!\[\matrix a,-N\\\frac {1} {2}+\frac {a} {2}-\frac {N} {2}\endmatrix; \frac
{1} {2}\]=\cases 0&\text {if $N$ is an odd nonnegative integer,}\\
\frac {\(\frac {1} {2}\)_{N/2}} {\(\frac {1} {2}-\frac {a}
{2}\)_{N/2}}&\text {if $N$ is an even nonnegative integer.}\endcases
\tag\BCi$$
As a result, in this case, the expression (\BCh) vanishes, whence
(\BCg) with $\ep=2$ is satisfied, and thus relation (\BCa)
specialized to the $(2n+2)$-nd column.

If $\ep=1$, the $_2F_1$-series in (\BCh) cannot be directly evaluated
by means of Gau{\ss}' formula. However, we may in a first stage apply
the contiguous relation
$$
{} _2 F _1 \!\left [ \matrix { a,b}\\ { c}\endmatrix ; {\displaystyle
   z}\right ]  = {} _2 F _1 \!\left [ \matrix { a,b - 1}\\ {
    c}\endmatrix ; {\displaystyle z}\right ]  + 
   {{az }\over
    {c}}\,
   {} _2 F _1 \!\left [ \matrix { a+1,b}\\ { c+1}\endmatrix ;
        {\displaystyle z}\right ]
$$
to transform (\BCh) into
$$\multline 
-\frac {(n-s-k)}2\cdot\frac {(-n+s+k+\frac {1} {2})_{2n-2s-2k-1}} 
{(2n-2s-2k)!}\\
-{\frac {   
       ({ \textstyle - n  + k + s+ {\frac 3 2}}) _{
         2 n - 2 s-2k-2}   }
   {8\,({2 n - 2 s-2k-2 })! }}
\Bigg({} _{2} F _{1} \!\left [ \matrix {1, - 2 n + 2 k +
        2 s+ 1}\\ {  -n + k + s+{\frac 3 2}}\endmatrix ;
        {\displaystyle {\frac 1 2}}\right ]\kern3cm\\
-
\frac {1} {2n - 2k - 2s- 3}
{} _{2} F _{1} \!\left [ \matrix {2, - 2 n + 2 k +
        2 s+ 2}\\ {  -n + k + s+{\frac 5 2}}\endmatrix ;
        {\displaystyle {\frac 1 2}}\right ]
\Bigg).
\endmultline$$
Both $_2F_1$-series in the last expression can now be evaluated by means
of Gau{\ss}' formula (\BCi). The first series simply vanishes, while
the second series evaluates to a non-zero expression. If this is
substituted, after some simplification we obtain
$$
-\frac {(-n+s+k+\frac {1} {2})_{2n-2s-2k-1}} 
{4\,(2n-2s-2k-1)!}
-{\frac {({ \textstyle -n + k + s+{\frac 3 2}}) _{2 n - 2 s-2k-2} } 
   {8\,({2 n - 2 s-2k-2})! }}=0.
$$
This shows that for $\ep=1$ the expression (\BCh) vanishes as well, whence
(\BCg) with $\ep=1$ is satisfied, and thus also relation (\BCa)
specialized to the $(2n+1)$-st column.

The verification of (\BCb) is completely analogous to that of (\BCe)
and is left to the reader.

\smallskip
{\smc Step 4. \it $\Pf M_n(x)$ is a polynomial in $x$ of degree
at most $2n^2+n-4k-3$}. 
By (\AB), $Q_{i,j}$ is a polynomial in $x$
of degree $i+j$. On the other hand, the degree of $H_{i,1}$ in $x$ is
clearly $i-2k-1$, while the degree of $H_{i,2}$ is $i-2k-2$. It
follows that, in the defining expansion of the determinant
$\det M_n(x)$, each nonzero term has degree 
$$
\sum _{i=1} ^{2n}i+
\sum _{j=1} ^{2n}j-2(2k+1)-2(2k+2)=4n^2+2n-8k-6.
$$
The Pfaffian being the square root of the determinant (cf.\ (\PF)), 
the claim follows.

\smallskip
{\smc Step 5. \it Evaluation of $P_n(x)$ at $x=0,-1,\dots,-n+k+1$}.
The polynomial $P_n(x)$ is defined by means of (\BD).
So, what we
would like to do is to set $x=-s$ in (\BD), $s$ being one of
$0,1,\dots,n-k-1$, evaluate $\Pf M_n(-s)$, divide both sides of (\BD) by
the products on the right-hand side of (\BD), and get the evaluation
of $P_n(x)$ at $x=-s$. However, the first
product on the right-hand side of
(\BD) unfortunately {\it is zero} for $x=-s$, $1\le s\le
n-k-1$. 
%Christian:
(It is not zero for $s=0$.)
Therefore we have to find a way around this difficulty.

%Mihai
%I agree, but shouldn't we say for completeness what the left hand side of 
%(5.3) evaluates to for s=0 and why?
%Christian: I have inserted a parenthetical remark above, and another remark
%  below.
%  I would not do more, because the only difference in case s=0 is
%  that this "preprocessing step" is not needed, but everything else
%  then is uniform in s, including s=0.
Fix an $s$ with $1\le s\le n-k-1$. Before setting $x=-s$ in (\BD), we
have to cancel $(x+s)^{s}$ (see (\BAa)) on both sides of (\BD). That is,
we should write (\BD) in the form
$$\multline
P_n(x)=
-\frac {1} {(x+s)^s}\Pf M_n(x)\\
\times
\underset \ell\ne s\to{\prod _{\ell=1} ^{n-k-1}}(x+\ell)^{-\ell}
\prod _{\ell=n-k+1} ^{n}(x+\ell)^{-\ell+1}
\prod _{\ell=n+1} ^{2n-1}(x+\ell)^{-2n+\ell-\chi(\ell>n+k)+1}\\
\times
{\prod _{\ell=1} ^{n}}(x+\ell+\tfrac {1} {2})_{2n-2\ell}^{-1},
\endmultline\tag\BF$$
and subsequently specialize $x=-s$. However, in order to be able to
perform this step, we need to evaluate
$$-\left(\frac {1} {(x+s)^{s}}\Pf M_n(x)\right)\bigg\vert_{x=-s}.$$
In order to accomplish this, we apply Lemma~\TCc\ with $N=2n+2$, $a=2n-2s$, 
$b=2n$, and $A=M_n(x)$. 
Indeed, $(x+s)$ is a factor of each entry in the $i$-th row in matrix
$M_n(x)$, for $i=2n-2s+1,2n-2s+2,\dots,2n$. We obtain
$$-\left(\frac {1} {(x+s)^{s}}\Pf M_n(x)\right)\bigg\vert_{x=-s}
=-\Pf(\widetilde Q)\,\Pf(S),\tag\BGa$$
where 
$$\widetilde Q=\pmatrix \overline Q&\overline H\\-\overline H^t&0\endpmatrix,
\tag\BG$$
with $\overline Q$ being given by
$$\overline Q=\left(Q_{i,j}\big\vert_{x=-s}\right)_{1\le i,j\le 2n-2s}$$
and $\overline H$ by
$$\overline H=\left(H_{i,j}\big\vert_{x=-s}\right)_{1\le i\le 2n-2s,\
1\le j\le 2},$$
and where
$$S=\left(
\left(\frac {1}
{x+s}Q_{i+2n-2s,j+2n-2s}\right)\bigg\vert_{x=-s}\right)_{1\le i,j\le 2s}.$$
%Christian:
We point out that (\BGa) also holds for $s=0$ once we interpret the
Pfaffian of an empty matrix (namely the Pfaffian of $S$) as $1$.
In particular, under that convention,
the arguments below can be used for $0\le s\le n-k-1$,
that is, {\it including} $s=0$.

We must now compute $\Pf(\widetilde Q)$ and $\Pf(S)$.
We start with the evaluation of $\Pf(S)$.
It follows from (\AB) that the $(i,j)$-entry of $S$ is given by
$$\multline 
S_{i,j}=
\sum _{l=0} ^{i+2n-2s-1}(-1)^{l+j+1}
\frac {j-i} {i+2n-2s}\binom
{j+2n-2s-1}{i+2n-2s-l-1}\binom {l+j+2n-2s}l\\
\cdot
\frac {(2n-2s+1)!\,(l+j-1)!}
{(l+j+2n-2s+1)!}.
\endmultline$$
If we write this using hypergeometric notation, we obtain the
alternative expression
$$S_{i,j}=
{\frac {{{\left( -1 \right) }^{j +1}} 
        ({ \textstyle j-i}) _{i + 2 n - 2 s} } 
    {(2 n - 2 s+j+1)! \,
      ({ \textstyle j}) _{i - j + 2 n - 2 s+1} }}
{} _{3} F _{2} \!\left [ \matrix { 1 - i - 2 n + 2 s, 1 + j + 2 n -
       2 s, j}\\ { 1 - i + j, 2 + j + 2 n - 2 s}\endmatrix ; {\displaystyle
       1}\right ].
$$
Rewrite this expression as the limit
$$S_{i,j}=\lim_{\ep\to0}
{\frac {{{\left( -1 \right) }^{j +1}} 
        ({ \textstyle j-i}) _{i + 2 n - 2 s} } 
    {(2 n - 2 s+j+1)! \,
      ({ \textstyle j}) _{i - j + 2 n - 2 s+1} }}
{} _{3} F _{2} \!\left [ \matrix { 1 - i - 2 n + 2 s, 1 + j + 2 n -
       2 s, j}\\ { 1 - i + j, 2+\ep + j + 2 n - 2 s}\endmatrix ; {\displaystyle
       1}\right ].
$$
Now we apply one of Thomae's
$_3F_2$-transformation formulas (cf\. \cite{\BailAA, Ex.~7, p.~98})
$$
{} _{3} F _{2} \!\left [ \matrix { a, b, c}\\ { d, e}\endmatrix ;
   {\displaystyle 1}\right ]  =
   \frac {\Ga( e)\,\Ga( d + e -a - b - c )} {\Ga(e -a)\,\Ga(  d +
    e-b - c)}
  {} _{3} F _{2} \!\left [ \matrix { a, -b + d, -c + d}\\ { d, -b - c + d +
    e}\endmatrix ; {\displaystyle 1}\right ]  .
$$
Thus, we obtain
$$\multline 
S_{i,j}=\lim_{\ep\to0}
{\frac {{{\left( -1 \right) }^{j + 1}} 
      \Gamma({ \textstyle 2 n - 2 s+\ep+1}) \,
      \Gamma({ \textstyle 2 n - 2 s+j+\ep+2}) \,
      ({ \textstyle j-i}) _{i + 2 n - 2 s} } 
    {\Gamma({ \textstyle \ep - i+2}) \,
      \Gamma({ \textstyle 4 n - 4 s+i+j+\ep+1}) \,
      (2 n - 2 s+j+1)! \,
      ({ \textstyle j}) _{i - j + 2 n - 2 s+1} }}\\
\times
      {} _{3} F _{2} \!\left [ \matrix { 1 - i - 2 n + 2 s, -i - 2 n +
       2 s, 1 - i}\\ { 1 - i + j, 2 + \ep - i}\endmatrix ;
       {\displaystyle 1}\right ],
\endmultline$$
or, in usual sum notation,
$$\multline 
S_{i,j}=\lim_{\ep\to0}
\sum _{l=0} ^{i-1} {\frac {{{\left( -1 \right) }^{j +1
  }}\,( j-i ) \,\Gamma({ \textstyle 2 n -
  2 s+\ep+1}) \,\Gamma({ \textstyle 2 n - 2 s+j+\ep+2}) } {\Gamma({
  \textstyle l- i + \ep+2}) \,\Gamma({ \textstyle 
   4 n - 4 s+i+j+\ep+1}) }}\\
\cdot
\frac {({
  \textstyle 1 - i}) _{l} \,({ \textstyle l - i + j + 1}) _{ 2 n
  - 2 s+i - l-1} \,({ \textstyle 2 n - 2 s+i - l}) _{l} } 
 {l! \,(2 n - 2 s+j+1)!
   \,({ \textstyle j}) _{2 n - 2 s+i-j-l+1} }.
\endmultline$$
Because of the term $\Gamma({ l- i + \ep+2})$ in the denominator, in
the limit only the summand for $l=i-1$ does not vanish. After
simplification, this leads to
$$S_{i,j}={\frac {{{\left( -1 \right) }^{i + j }}\,
      ( j-i ) \,(2n-2s+i-1)!\,(2n-2s+j-1)! }
      {
(4n-4s+i+j)!\,(2n-2s+1)! }}.
$$
We must evaluate the Pfaffian
$$\underset {1\le i,j\le 2s}\to\Pf(S_{i,j}).$$
By factoring some terms out of rows and columns, we see that
$$\multline
\underset {1\le i,j\le 2s}\to\Pf(S_{i,j})=(-1)^s
{(2n-2s+1)!}^{-s}\\
\times\bigg(\prod _{i=1} ^{2s} {
(2n-2s+i-1)!} 
\bigg)
\underset {1\le i,j\le 2s}\to\Pf\(\frac {j-i} {(4n-4s+i+j)!}\).
\endmultline$$
This Pfaffian can be evaluated in closed form by Corollary~\TCb\ in
the next section.
The result is that
$$
\Pf(S)=(-1)^s
{(2n-2s+1)!}^{-s}
\bigg(\prod _{i=1} ^{2s}{
(2n-2s+i-1)!} 
\bigg)
\bigg(
\prod _{i=0} ^{s-1}\frac {(2i+1)!} {(4n-2s+2i+1)!}\bigg).
\tag\BH$$

\midinsert
\vskip10pt
\vbox{\noindent
\centertexdraw{
\drawdim truecm \setunitscale.8
\linewd.15 
\RhombusA \RhombusA \RhombusA \RhombusA \RhombusA
\RhombusA \RhombusA \RhombusA \RhombusA \RhombusA 
\move(.866025 .5)
\RhombusA \RhombusA \RhombusA \RhombusA \RhombusA 
\RhombusA \RhombusA \RhombusA \RhombusA 
\move(1.73205 1)
\RhombusA \RhombusA 
\move(2.59808 1.5)
\RhombusA \RhombusA 
\move(3.464 2)
\RhombusA \RhombusA 
\move(4.33013 2.5)
\RhombusA \RhombusA 
\move(5.19615 3)
\RhombusA \RhombusA 
\move(6.06218 3.5)
\RhombusA \RhombusA 
\move(6.9282 4)
\RhombusA \RhombusA 
\move(7.79423 4.5)
\RhombusA 
\move(5.19615 -1)
\RhombusA \RhombusA \RhombusA \RhombusA 
\move(5.19615 1)
\RhombusB \RhombusB
\move(6.06218 1.5)
\RhombusA \RhombusB \RhombusB \RhombusA \RhombusA 
\move(6.9282 2)
\RhombusA \RhombusB \RhombusA 
\move(7.79423 2.5)
\RhombusA 
\move(7.79423 0.5)
\RhombusC 
\move(6.06218 0.5)
\RhombusC 
\move(10 0)
\bsegment
\drawdim truecm \linewd.15
\RhombusA \RhombusA \RhombusA \RhombusA \RhombusA
\RhombusA \RhombusA \RhombusA \RhombusA \RhombusA 
\move(.866025 .5)
\RhombusA \RhombusA \RhombusA \RhombusA \RhombusA 
\RhombusA \RhombusA \RhombusA \RhombusA 
\move(1.73205 1)
\RhombusA \RhombusA 
\move(2.59808 1.5)
\RhombusA \RhombusA 
\move(3.464 2)
\RhombusA \RhombusA 
\move(4.33013 2.5)
\RhombusA \RhombusA 
\move(5.19615 3)
\RhombusA \RhombusA 
\move(6.06218 3.5)
\RhombusA \RhombusA 
\move(6.9282 4)
\RhombusA \RhombusA 
\move(7.79423 4.5)
\RhombusA 
\linewd.02
\move(5.19615 0)
\rdreieck \rhombus \rhombus \rhombus 
\move(5.19615 1)
\rdreieck \rhombus \rhombus \rhombus 
\move(5.19615 1)
\rhombus \rhombus \rhombus 
\move(6.06218 1.5)
\rhombus \rhombus
\move(6.9282 2)
\rhombus 
\move(8.66025 -2)
\vdSchritt
\move(8.66025 -1)
\vdSchritt
\move(8.66025 0)
\vdSchritt
\move(8.66025 1)
\vdSchritt
\move(8.66025 2)
\vdSchritt
\move(8.66025 3)
\vdSchritt
\esegment
}
\centerline{\eightpoint a. A lozenge tiling for the degenerate region\quad \quad 
\quad \quad 
b. Forced lozenges in case $x=0$}
\vskip8pt
\centerline{\eightpoint Figure \FD}
}
\vskip10pt
\endinsert

\medskip
We finally turn to the evaluation of $\det(\widetilde Q)$. 
If we compare (\BG) with (\AAAA), then we see that $\widetilde
Q=M_{n-s}(0)$. Hence, using Lemma~\TAE\ with $n$ replaced by $n-s$ and
with $x=0$, we see that 
$-\Pf(\widetilde Q)$ 
is equal to $\M(F_{n-s,0}\setminus\triangleleft_2(k))$. 
(The reader should recall the definitions of the region $F_{n,x}$ 
and of the triangular hole $\triangleleft_2(k)$ given in the
introduction, see again Figure~\FF.) 
Figure~\FD.a shows a typical example 
where $n-s=5$ and $k=2$. Since this region is degenerate, there are 
many forced lozenges, see Figure~\FD.b. The enumeration problem
therefore reduces to the problem of determining the 
number of symmetric lozenge tilings 
of a hexagon with side lengths $2k,2k,2,2k,2k,2$. This number is
given by formula ({\SPP}) with $n=k$ and $x=1$. If we substitute this in
(\BGa), together with the evaluation (\BH), then,
after some manipulation, we obtain 
$$\multline 
-\left(\frac {1} {(x+s)^{s}}\Pf M_n(x)\right)\bigg\vert_{x=-s}=(-1)^s
    \binom {4k+1}{2k}
\frac {(2s)!} {(2n-2s+1)!^{s}\,2^{s}\,s!}\\
\times
\bigg(\prod _{i=1} ^{2s}(2n-2s+i-1)!\bigg)
\bigg(\prod _{i=0} ^{s-1}\frac {(2i)!} {(4n-2s+2i+1)!}\bigg).
\endmultline$$
Hence, by inserting this in (\BF), we have
$$\multline
P_n(-s)=(-1)^s
    \binom {4k+1}{2k}
\frac {(2s)!} {(2n-2s+1)!^{s}\,2^{s}\,s!}\\
\times
\bigg(\prod _{i=1} ^{2s}(2n-2s+i-1)!\bigg)
\bigg(\prod _{i=0} ^{s-1}\frac {(2i)!} {(4n-2s+2i+1)!}\bigg)\\
\times
\underset \ell\ne s\to{\prod _{\ell=1} ^{n-k-1}}(-s+\ell)^{-\ell}
\prod _{\ell=n-k+1} ^{n}(-s+\ell)^{-\ell+1}
\prod _{\ell=n+1} ^{2n-1}(-s+\ell)^{-2n+\ell-\chi(\ell>n+k)+1}\\
\times
{\prod _{\ell=1} ^{n}}(-s+\ell+\tfrac {1} {2})_{2n-2\ell}^{-1}.
\endmultline\tag\BI$$

This completes the proof of Theorem~{\TB}.

\head 6. An auxiliary determinant evaluation,
and an auxiliary Pfaffian factorization \endhead
Mehta and Wang proved the following determinant
evaluation in \cite{\MeWaAA}. (There is a typo in the formula stated in
\cite{\MeWaAA, Eq.~(7)} in that the binomial coefficient $\binom nk$
is missing there.)

\proclaim{Theorem \TAF}
For all real numbers $a,b$ and positive integers $n$, we have
$$\multline 
\det_{0\le i,j\le n-1}\big((a+j-i)\,\Gamma(b+i+j)\big)\\=
     \bigg(\prod _{i=0} ^{n-1}i!\,\Gamma(b+i)\bigg)
\sum _{k=0} ^{n}(-1)^k \binom nk \big((b-a)/2\big)_k\,
 \big((b+a)/2\big)_{n-k},
\endmultline\tag\CAa$$
as long as the arguments occurring in the gamma functions avoid their
singularities. 
\endproclaim

The sum on the right-hand side of (\CAa) can be alternatively
expressed as the coefficient of $z^n/n!$ in
$$(1+z)^{(a-b)/2}(1-z)^{(-a-b)/2}.$$
Therefore,
in the case $a=0$ we obtain the following simpler determinant
evaluation.

\proclaim{Corollary \TAG}
For all real numbers $b$ and positive integers $n$, we have
$$
\det_{0\le i,j\le n-1}\big((j-i)\,\Gamma(b+i+j)\big)=
   \chi(n\text { \rm is even})
   \bigg(\prod _{i=0} ^{n-1}i!\,\Gamma(b+i)\bigg)
\frac {n!\,(b/2)_{n/2}} {(n/2)!},
$$
as long as the arguments occurring in the gamma functions avoid their
singularities. Here, as before, $\chi(\Cal A)=1$ if $\Cal A$ is
true and $\chi(\Cal A)=0$ otherwise.
\endproclaim

One can obtain the following slightly
(but, for our purposes, essentially) stronger
statement. It is stated as Eq.~(4) in \cite{\MeWaAA}, with the
argument how to obtain it hinted at at the bottom of page~231 of
\cite{\MeWaAA}. Since, from there, it is not completely obvious how to
actually complete the argument, we provide a proof.

\proclaim{Proposition \TCa}
For all real numbers $b$ and positive even integers $n$, we have
$$
\underset {0\le i,j\le n-1}\to{\Pf}\big((j-i)\,\Gamma(b+i+j)\big)=
     \prod _{i=0} ^{\frac n2-1}(2i+1)!\,\Gamma(b+2i+1),
\tag\CAb$$
as long as the arguments occurring in the gamma functions avoid their
singularities. 
\endproclaim

\demo{Proof}
Since the Pfaffian of a skew-symmetric matrix equals the square root
of its determinant (cf.\ (\PF)), the formula given by Theorem~\TAG\ yields, after a
little manipulation, that
$$
\underset {0\le i,j\le n-1}\to{\Pf}\big((j-i)\,\Gamma(b+i+j)\big)=
    \ep \prod _{i=0} ^{\frac n2-1}(2i+1)!\,\Gamma(b+2i+1),
\tag\CAc$$
where $\ep=+1$ or $\ep=-1$. In order to determine the sign $\ep$, we
argue by induction on (even) $n$. Let us suppose that we have already proved 
(\CAb) up to $n-2$. We now multiply both sides of (\CAc) by $b+1$ and
then let $b$ tend to $-1$. Thus, on the right-hand side we obtain the
expression 
$$
\ep\bigg(\prod _{i=0} ^{\frac n2-1}(2i+1)!\bigg)\bigg(
\prod _{i=1} ^{\frac n2-1}\Gamma(2i)\bigg).
\tag\CAe$$
On the other hand, by the definition  of the Pfaffian, on the
left-hand side we obtain
$$
\sum _{\pi\in \Cal M[0,\dots, n-1]} ^{}
\sgn\pi
\lim_{b\to-1}\Bigg((b+1)
\underset i,j\text{ matched in }\pi\to{\prod _{i<j} ^{}}
(j-i)\,\Gamma(b+i+j)\Bigg)
\tag\CAd$$
(with the obvious meaning of $\Cal M[0,\dots, n-1]$; cf\. the sentence
containing (\PPF)).
In this sum, matchings $\pi$ for which all matched pairs $i,j$
satisfy $i+j>1$ do not contribute, because the corresponding summands
vanish. However, there is only one possible pair $i,j$ with $0\le i<j$
for which $i+j\le1$, namely $(i,j)=(0,1)$. 
Therefore, the sum in (\CAd) reduces to
$$\align
\sum _{\pi'\in \Cal M[2,\dots, n-1]} ^{}&
\sgn\pi'
\Big(\lim_{b\to-1}(b+1)(1-0)\,\Gamma(b+1)\Big)
\underset i,j\text{ matched in }\pi'\to{\prod _{i<j} ^{}}
(j-i)\,\Gamma(i+j-1)\\
&=
\underset {2\le i,j\le n-1}\to{\Pf}\big((j-i)\,\Gamma(i+j-1)\big)\\
&=
\underset {0\le i,j\le n-3}\to{\Pf}\big((j-i)\,\Gamma(i+j+3)\big),
\endalign$$
where the next-to-last equality holds by the definition (\PPF) of the
Pfaffian. Now we can use the induction hypothesis to evaluate the last
Pfaffian. Comparison with (\CAe) yields that $\ep=+1$.\quad \quad \qed
\enddemo

By using the reflection formula (cf\. \cite{\AnARAA, Theorem~1.2.1})
$$\Gamma(x)\,\Gamma(1-x)=\frac {\pi} {\sin \pi x}$$
for the gamma function, and the substitutions $i\to n-i-1$ and
$j\to n-j-1$, it is not difficult to see that Proposition~\TCa\
is equivalent to the following.

%\proclaim{Lemma \TCa} 
%For all positive integers $n$, we have
%$$\multline 
%\det_{0\le i,j\le n-1}\left(\frac {a+i-j} {\Gamma(b+i+j)}\right)=
%   (-1)^{\binom {n+1}2} 
%   \Bigg(\prod _{i=0} ^{n-1}\frac {i!} {\Gamma(b+n+i-1)}\Bigg)\\
%\times
%\sum _{k=0} ^{n}(-1)^k \binom nk \big(\tfrac {1} {2}(3+a-b-2n)\big)_k\,
% \big(\tfrac {1} {2}(3-a-b-2n)\big)_{n-k}.
%\endmultline\tag\CAa$$
%\endproclaim
%
%Since the sum on the right-hand side of (\CAa) can be alternatively
%expressed as the coefficient of $z^n/n!$ in
%$$(1+z)^{(-3-a+b+2n)/2}(1-z)^{(-3+a+b+2n)/2},$$
%in the case $a=0$ we obtain the following simpler result.

\proclaim{Corollary \TCb}
For all positive even integers $n$, we have
$$
\underset {0\le i,j\le n-1}\to{\Pf}\left(\frac {j-i} {\Gamma(b+i+j)}\right)=
   \prod _{i=0} ^{\frac {n} {2}-1}\frac {(2i+1)!} {\Gamma(b+n+2i-1)}.
$$
\endproclaim

We close this section by proving a factorization of a 
certain specialization of a Pfaffian that we need
in Step~5 in Section~5.

\proclaim{Lemma \TCc}
Let $N,a,b$ be positive integers with $a< b\le N$, where $N$ and
$b-a$ are even.
Let $A=(A_{i,j})_{1\le i,j\le N}$ be a skew-symmetric matrix
with the following properties:

\roster 
\item The entries of $A$ are polynomials in $x$.
\item The entries in rows $a+1,a+2,\dots,b$ {\rm(}and, hence, also in
the corresponding columns{\rm)} are divisible by $x+s$.
\endroster

Then
$$\left(\frac {1} {(x+s)^{(b-a)/2}}\Pf A\right)\bigg\vert_{x=-s}=\Pf \widetilde
A\cdot\Pf S,
\tag\CAg$$
where $\widetilde A$ is the matrix which arises from $A$ by deleting
rows and columns $a+1,a+2,\dots,b$ and subsequently specializing
$x=-s$, and 
$$S=\left(
\left(\frac {1} {x+s}A_{i,j}\right)\bigg\vert_{x=-s}\right)
_{a+1\le i,j\le b}.$$
\endproclaim

\demo{Proof}
By the definition (\PPF) of the Pfaffian, we have
$$\left(\frac {1} {(x+s)^{(b-a)/2}}\Pf A\right)\bigg\vert_{x=-s}
=
\Bigg(\frac {1} {(x+s)^{(b-a)/2}}\sum _{\pi\in\Cal M[1,\dots, N]} ^{}\sgn\pi 
\underset i,j\text{ matched in }\pi\to{\prod _{i<j} ^{}}
A_{i,j}\Bigg)\Bigg\vert_{x=-s}.
$$
Let $\Cal M_1$ denote the subset of $\Cal M[1,\dots, N]$ consisting of
those matchings that pair all the elements from $\{a+1,a+2,\dots,b\}$
among themselves (and, hence, all the elements of the complement
$\{1,2,\dots,a,b+1,b+2,\dots,N\}$ among themselves). Let
$\Cal M_2$ be the complement $\Cal M[1,\dots, N]\backslash \Cal M_1$. 
Then
$$\multline
\left(\frac {1} {(x+s)^{(b-a)/2}}\Pf A\right)\bigg\vert_{x=-s}
=
\Bigg(\frac {1} {(x+s)^{(b-a)/2}}\sum _{\pi\in\Cal M_1} ^{}\sgn\pi 
\underset i,j\text{ matched in }\pi\to{\prod _{i<j} ^{}}
A_{i,j}\Bigg)\Bigg\vert_{x=-s}\\
+
\Bigg(\frac {1} {(x+s)^{(b-a)/2}}\sum _{\pi\in\Cal M_2} ^{}\sgn\pi 
\underset i,j\text{ matched in }\pi\to{\prod _{i<j} ^{}}
A_{i,j}\Bigg)\Bigg\vert_{x=-s}.
\endmultline
\tag\CAf$$
Each term in the sum in the second line of (\CAf) vanishes, since
the product contains more than $(b-a)/2$ factors $A_{i,j}$ that are
divisible by $x+s$. On the other hand, every matching $\pi$ in $\Cal M_1$
is the disjoint union of a matching $\pi'\in 
M[1,2,\dots, a,b+1,b+2,\dots,N]$ and a matching 
$\pi''\in M[a+1,a+2,\dots,b]$. If we also use the
simple fact that $\sgn \pi=\sgn\pi'\cdot\sgn\pi''$
(as there are no crossings between paired elements of $\pi'$
and paired elements of $\pi''$),
then we obtain
$$\align
\bigg(\frac {1} {(x+s)^{(b-a)/2}}&\Pf A\bigg)\bigg\vert_{x=-s}\\
&=
\left(\frac {1} {(x+s)^{(b-a)/2}}
\underset \pi''\in  M[a+1,\dots,b]\to{\sum _{\pi'\in\Cal M[1,\dots, a,b+1,\dots,
N]} ^{}}\sgn\pi'\cdot\sgn\pi''\right. \\
&\kern2cm
\left.\left.\cdot
\bigg(\underset i,j\text{ matched in }\pi'\to{\prod _{i<j} ^{}}
A_{i,j}\bigg)
\bigg(
\underset i,j\text{ matched in }\pi''\to{\prod _{i<j} ^{}}
A_{i,j}\bigg)\right)\right\vert_{x=-s}\\
&=
{\sum _{\pi'\in\Cal M[1,\dots, a,b+1,\dots,
N]} ^{}}\sgn\pi'
\underset i,j\text{ matched in }\pi'\to{\prod _{i<j} ^{}}
A_{i,j} \Big\vert_{x=-s}\\
&\kern2cm
\cdot
{\sum _{ \pi''\in  M[a+1,\dots,b]} ^{}}\sgn\pi''
\underset i,j\text{ matched in }\pi''\to{\prod _{i<j} ^{}}
\left(\frac {1} {x+s}A_{i,j}\right)\bigg\vert_{x=-s}.
\endalign
$$
By the definition (\PPF) of the Pfaffian,
the last expression is exactly the right-hand side of (\CAg).\quad
\quad \qed
\enddemo

\head 7. Proofs of Theorems {\TAA} and {\TAB} \endhead

In our proofs we make use of the following lemmas.

\proclaim{Lemma \TC}Let $\be$ be a real number with either $\be>0$ or
$\be<-1$. Then, for all sequences $(\be_n)_{n\ge1}$ with $\be_n\to\be$
as $n\to\infty$, we have
$$
\lim_{n\to\infty}\frac {1} {\sqrt n}\,\,
{} _{5} F _{4} \!\left [ \matrix { -2n,\frac {1}
{2},-n+k+1,-n+k+1,\be_n n}\\ 
{ -2n+\frac {1} {2},-n-k,-n-k,\be_n n+1}\endmatrix ;
   {\displaystyle 1}\right ]=
\int _{0} ^{1}
\frac {\sqrt 2\,(1-\al)^{4k+2}} {(1+\frac {\al} {\be})
\,\sqrt{\pi\al(2-\al)}}\,d\al,
\tag\CA$$
where, on the left-hand side, we used again the standard notation
{\rm(\HYP)} for hypergeometric series.
\endproclaim
\demo{Proof} 
We write the $_5F_4$-series in (\CA) explicitly as a sum over $l$:
$$
\sum _{l=0} ^{n-k-1}\frac {\Ga(2n+1)\,\Gamma(l+\frac {1} {2})\,
\Ga(2n-l+\frac {1} {2})\,\Ga(n-k)^2\,\Ga(n+k-l+1)^2} 
{\Ga(2n-l+1)\,\Ga(\frac {1} {2})\,\Ga(l+1)\,\Ga(2n+\frac 12)\,\Ga(n-k-l)^2\,\Ga(n+k+1)^2}
\frac {\be_n n} {(\be_n n+l)}.\tag\CB$$
Let us denote the summand in this sum by $F(n,l)$. We have
$$\multline
\frac {\partial} {\partial l}F(n,l)=F(n,l)
\bigg(\psi(l+\tfrac12)-\psi(l+1)+
        \psi(2n-l+1)-\psi(2n-l+\tfrac12)\\
+2\psi(n-k-l)-2\psi(n+k-l+1)-
\frac {1} {\be_n n+l}\bigg),
\endmultline$$
where $\psi(x):=(\frac {d} {dx}\Ga(x))/\Ga(x)$ is the digamma function. Since
$\psi(x)$ is a monotone increasing, concave function for $x>0$
(this follows e.g\. from \cite{\AnARAA, Eq.~(1.2.14)}), we
have 
$$\psi(l+1)-\psi(l+\tfrac12)\ge\psi(2n-l+1)-\psi(2n-l+\tfrac12)$$
for $0\le l\le n$. Moreover, because of the equality
$\psi(x+1)=\psi(x)+\frac {1}
{x}$ (cf\. \cite{\AnARAA, Eq.~(1.2.15) with $n=1$}), and since either
$\be>0$ or $\be<-1$, for large enough $n$ we have
$$\psi(n+k-l+1)\ge\psi(n-k-l)+\frac {1} {n+k-l}>\psi(n-k-l)-\frac {1}
{\be_n n+l}.$$
Altogether, this implies that $\frac {\partial} {\partial l}F(n,l)<0$
for $0\le l\le n-k-1$, that is, for fixed large enough $n$, the summand $F(n,l)$
is monotone decreasing as a function in $l$. In particular, 
for $0\le l\le n-k-1$ we have
$$0<F(n,l)\le F(n,0)=1.\tag\CC$$
The sum (\CB) may therefore be approximated by an integral:
$$\align \sum _{l=0} ^{n-k-1}F(n,l)&=
\sum _{l=0} ^{\fl{\log n}-1}F(n,l)+\sum _{l=\fl{\log n}}
^{n-k-\fl{\log n}-1}F(n,l)+\sum _{l=n-k-\fl{\log n}} ^{n-k-1}F(n,l)\\
%Christian: no "funny" term
%&=O(\log n)+\int _{ \fl{\log n}-1} ^{n-k-\fl{\log n}-1}F(n,l)\,d\,l+
%  O\big(F(n,\fl{\log n}-1)\big)
%\\
&=O(\log n)+\int _{ \fl{\log n}-1} ^{n-k-\fl{\log n}-1}F(n,l)\,d\,l.
\tag\CD
\endalign$$

%Mihai: Last term in middle line above should be O(\log n), right?
%Christian: The "funny" term $O\big(F(n,\fl{\log n}-1)$ is supposed to
% bound the error when you replace the sum over l in the center by
% the integral. It is easy to see that the monotonicity of the summand
% implies that this error is bounded above by the value of the function
% at the lower limit of integration, that is, at $l=\fl{\log n}-1$,
% whence the funny expression. The expression itself is bounded by $O(1)$
% because of (7.8).
% I do not know whether it makes sense to add any explanations, whether
% one should not mention this at all (and simply delete the funny expression).
%Mihai
%I see, you're right. It may be clearer to drop the ``funny'' term altogether.
%Christian: okay.

The next step is to apply Stirling's approximation
$$\log\Ga(z)=\(z-\frac {1} {2}\)\log(z)-z+\frac {1} {2}\log(2\pi)
+O\(\frac {1} {z}\)
\tag\CE$$
for the gamma function, in the form
$$\align 
\log\Ga(an+bl+c)&=\(a n+b l+c-\frac {1} {2}\)\(\log\(a+b \tfrac
{l\vphantom{1}} {n}\)
  +\log(n)+\log\(1+\tfrac {c} {an+bl}\)\)\\
&\kern1cm
        -(a n+b l+c)+
        \frac {1} {2}\log(2\pi)+O\(\frac {1} {an+bl}\)\\
&=\(a n+b l+c-\frac {1} {2}\)(\log(a+b \tfrac ln)+\log(n))\\
&\kern1cm
-(a n+b l)+
        \frac {1} {2}\log(2\pi)+O\(\frac {1} {an+bl}\),
\endalign$$
where $a,b,c$ are real numbers with $a\ge 0$.
If this is used in the defining expression for $F(n,l)$, then 
after cancellations we obtain
$$\align 
\log F(n,l)&=  \frac {1} {2}{{\log (2)}} + 
    \left( 4 k + 2 \right)  \log \(1 - {\frac l n}\) - 
   \frac {1} {2}{{\log \(2 - {\frac l n}\)}}\\
&\kern3cm
- \frac {1} {2}{{\log \( {\frac l n}\)}}
- \frac {1} {2}{{\log ( n)}}
- \frac {1} {2}{{\log (\pi)}} - \log \(1 + {\frac l {\be_n
n}}\)\\
&\kern3cm
+O\(\frac {1} {l}\)+O\(\frac {1} {n-l}\)+O\(\frac {1} {2n-l}\)+
O\(\frac {1} {n}\)\\
&=\log\(\frac {\sqrt 2\,(1-\frac {l} {n})^{4k+2}} {\sqrt n(1+\frac {l} {\be_n n})
\,\sqrt{\pi\frac {l} {n}(2-\frac {l} {n})}}\)
+O\(\frac {1} {\log n}\),
\endalign$$
as long as $\log n\le l\le n-k-\log n$. 
Substitution of this approximation in (\CD) yields
$$\multline \sum _{l=0} ^{n-k-1}F(n,l)\\=
\Bigg(\int _{ \fl{\log n}-1} ^{n-k-\fl{\log n}-1}
\frac {\sqrt 2\,(1-\frac {l} {n})^{4k+2}} {\sqrt n(1+\frac {l} {\be_n n})
\,\sqrt{\pi\frac {l} {n}(2-\frac {l} {n})}}\,d\,l\Bigg)
\(1+O\(\frac {1} {\log n}\)\)+O\({\log n}\),
\endmultline$$
or, after the substitution $l=\al n$,
$$\multline \sum _{l=0} ^{n-k-1}F(n,l)\\=
\sqrt n\Bigg(\int _{(\fl{\log n}-1)/n} ^{(n-k-\fl{\log n}-1)/n}
\frac {\sqrt 2\,(1-\al)^{4k+2}} {(1+\frac {1} {\be_n}\al)
\,\sqrt{\pi\al(2-\al)}}\,d\al\Bigg)
\(1+O\(\frac {1} {\log n}\)\)+O\({\log n}\).
\endmultline$$
The assertion of the lemma follows now immediately.\quad \quad \qed
\enddemo

With no extra work we can now get an exact formula for a
generalization $\omega_f(k;\xi)$ of the correlation 
$\omega_f(k)$ described in Section~2. For any real number $\xi>0$,
define $\omega_f(k;\xi)$ in analogy to  
({\Bmm}) by
$$
\omega_f(k;\xi):=\lim_{n\to\infty}\frac{\M(F_{n,\xi_nn}\setminus
\triangleleft_2(k))}{\M(F_{n,\xi_nn})},
\tag\X
$$
where $(\xi_n)_{n\ge1}$ is a suitable sequence of rational numbers
approaching $\xi$. (``Suitable" here means that we have to choose
$\xi_n$ in such a way that $\xi_nn$ is integral.) The number
$\omega_f(k;\xi)$ is the correlation of the triangular gap
$\triangleleft_2(k)$ with the free boundary, obtained  
when the large regions used in the definition are the left halfs of
hexagons that are not necessarily regular, but have their left
vertical side $\xi$ times as long as the two oblique sides. 
Note that,
by the results of \cite{\CLP}, we should expect distorted dimer
statistics around the gap for $\xi\neq1$. As Theorem~{\TF} shows, 
the distortion is quite radical: $\omega_f(k;\xi)$ turns out to decay
exponentially to 0 or blow up exponentially, 
according as $\xi>1$ or $\xi<1$; see also Remark~1.

\proclaim{Lemma \TD}
For any $\xi>0$ and $0\leq k\in\Z$, we have
$$\multline 
\om_f(k;\xi)=
\frac{1}{\pi}\,\binom {4k+1}{2k}\,\frac{1}{(1+\xi)^{4k+2}\sqrt{2+\xi}}\\
\times
\Bigg((\xi+2) \int _{0} ^{1}
\frac {(1-\al)^{4k+2}} {(1+\frac {\al} {\xi})
\,\sqrt{\al(2-\al)}}\,d\al -
\xi\int _{0} ^{1}
\frac {(1-\al)^{4k+2}} {(1-\frac {\al} {2+\xi})
\,\sqrt{\al(2-\al)}}\,d\al
\Bigg)\\
=\frac{1}{\pi}\,\binom {4k+1}{2k}\,\frac{1}{(1+\xi)^{4k+2}\sqrt{2+\xi}}
\int _{0} ^{1}
\frac {2\,(1-\al)^{4k+3}} {(1+\frac {\al} {\xi})
\,(1-\frac {\al} {2+\xi})\,\sqrt{\al(2-\al)}}\,d\al.
\endmultline\tag\CF$$
\endproclaim

\demo{Proof} 
By Theorem~{\TB} and formula~(\SPP), the ratio between
$\M(F_{n,x}\setminus\triangleleft_2(k))$ and 
$\M(F_{n,x})$
is, when written in hypergeometric notation,
$$\multline 
\binom {4k+1}{2k}\,
\frac {(n+k)!\,(\frac {1} {2})_{2n}} {(x+n-k)_{2k+1}\,(x+\frac {1} {2})_{2n}}
\frac {(x+1)_{n-k-1}\,(x+n+k+1)_{n-k-1}  }
  {(n-k-1)!^2\,(n+k+1)_{n-k}}\\
\times
\Bigg((x+2n)\, {} _{5} F _{4} \!\left [ \matrix { -2n,\frac {1}
{2},-n+k+1,-n+k+1,x}\\ 
{ -2n+\frac {1} {2},-n-k,-n-k,x+1}\endmatrix ;
   {\displaystyle 1}\right ]\\-
x\,{} _{5} F _{4} \!\left [ \matrix { -2n,\frac {1}
{2},-n+k+1,-n+k+1,-2n-x}\\ 
{ -2n+\frac {1} {2},-n-k,-n-k,-2 n-x+1}\endmatrix ;
   {\displaystyle 1}\right ]
\Bigg).
\endmultline$$
We now substitute $x=\xi_n n$ in this expression. Use of Lemma~\TC\
(which applies, as $\xi>0$),
together with Stirling's formula (\CE), yields the assertion.\quad \quad
\qed
\enddemo

\proclaim{Lemma \TE}
For any $\beta\neq0$ we have
$$
\int _{0} ^{1}
\frac {(1-\al)^{4k+2}} {(1+\frac {\al} {\beta})\,\sqrt{\al(2-\al)}}\,d\al 
\sim \sqrt{\frac{\pi}{8 k}},\ \ \ k\to\infty.\tag\CG
$$
\endproclaim

\demo{Proof} 
Let $I_{\beta}(k)$ be the integral on the left hand side of (\CG). The
asymptotics of $I_{\beta}(k)$ as $k\to\infty$ 
can be readily found using Laplace's method as presented for instance
in \cite{\OlveAA}. Conditions  
$(i)$--$(v)$ of \cite{\OlveAA, pp.~121--122} are readily checked. By
\cite{\OlveAA, Theorem~6.1, p.~125}, the 
large $z$ asymptotics of $\int_a^b e^{-z p(t)}q(t)\,dt$ is determined
by the quantities $\lambda$, $\mu$, $p_0$  
and $q_0$ in the series expansions
$$
p(t)-p(a)=p_0(t-a)^{\mu}+p_1(t-a)^{\mu+1}+\cdots
$$
and
$$
q(t)=q_0(t-a)^{\lambda}+q_1(t-a)^{\lambda+1}+\cdots.
$$
Namely, under the above assumptions one has
$$
e^{z p(a)} \int_a^b e^{-z p(t)}q(t)\,dt = 
\Gamma\left(\frac{\lambda}{\mu}\right)
\frac{q_0/(\mu p_0^{\lambda/\mu})}{z^{\lambda/\mu}}
+O\left(\frac{1}{z^{\lambda/\mu+1}}\right).\tag\CH
$$
In the case of $I_{\beta}(k)$ we have $p(t)=-\ln(1-t)$,
$q(t)=\frac{1}{(1-t/\beta)\sqrt{t(2-t)}}$, $a=0$, and $b=1$.  
These yield parameters $\lambda=1/2$, $\mu=1$, $p_0=1$, and
$q_0=1/\sqrt{2}$. In addition, $p(a)=0$. As in our case $z=4k+2$,
under these specializations ({\CH}) becomes ({\CG}). 
\quad \quad
\qed\enddemo

\proclaim{Theorem \TF}
As $k\to\infty$, the correlation $\om_f(k;\xi)$ is asymptotically
$$
\om_f(k;\xi)\sim \frac {1} {\pi(1+\xi)^2 \sqrt{\xi(2+\xi)}}\cdot 
\frac {1} {k}\(\frac {2}
{1+\xi}\)^{4k}.
$$
\endproclaim

\demo{Proof} 
Combine Lemmas~\TD\ and \TE\ with Stirling's approximation for 
the binomial coefficient $\binom {4k+1}{2k}$ in (\CF).\quad \quad
\qed
\enddemo

\medskip
\flushpar
{\it Proof of Theorem~{\TAA}.} Set $\xi=1$ in Theorem~{\TF}. \quad \quad \qed 

\medskip
\flushpar
{\it Proof of Theorem~{\TAB}.} Set
$$
D_k:=3 I_{1}(k)- I_{-3}(k),
$$
where $I_{\beta}(k)$ denotes the integral on the left hand side of
(\CG). Recalling that $\omega_f(k)$ is the $\xi=1$ specialization of
$\omega_f(k;\xi)$, we have by Lemma~{\TD} that  
$$
\omega_f(k+1)-\omega_f(k)=
\frac{1}{\pi}\frac{1}{2^{4k+2}\sqrt{3}}\binom{4k+1}{2k}
\left\{\left[\frac{(4k+3)(4k+5)}{4(2k+2)(2k+3)}-1\right]D_{k+1}+
(D_{k+1}-D_k)\right\},
$$
and thus
$$
\frac{\omega_f(k+1)-\omega_f(k)}{\omega_f(k)}=
\left[\frac{(4k+3)(4k+5)}{4(2k+2)(2k+3)}-1\right]
\frac{D_{k+1}}{D_{k}}+\frac{D_{k+1}-D_{k}}{D_{k}}.
\tag\CI
$$
By two applications of Lemma~{\TE} it follows that 
$$
D_k\sim\frac{\sqrt{\pi}}{\sqrt{2}k},\ \ \ k\to\infty.\tag\CJ
$$
Thus $D_{k+1}/D_k\to1$ as $k\to\infty$, and elementary arithmetics 
implies that the first term on the right hand side of ({\CI}) is
asymptotically $-1/(2k)$ as $k\to\infty$.  

To determine the asymptotics of the second term, write by Lemma~{\TD}
$$
D_{k+1}-D_k=3\left[I_{1}(k+1)-I_{1}(k)\right]-\left[I_{-3}(k+1)-I_{-3}(k)\right].
\tag\CK
$$
%If $\beta>0$ or $\beta<-1$ we have by Lemma~{\TD} that 
As $I_\beta(k)$ is the integral on the left hand side of (\CG), we have
$$
I_{\beta}(k+1)-I_{\beta}(k)=
\int_0^1\frac{(1-\alpha)^{4k+2}}{\left(1+\frac{\alpha}{\beta}\right)\sqrt{\alpha(2-\alpha)}}
\left[(1-\al)^4-1\right]d\alpha.\tag\CL
$$ 
The asymptotics of the integral in ({\CL}) follows by Laplace's
method, in the same manner as the proof of Lemma~{\TE}. In this case
$\lambda=3/2$, $\mu=1$, and equations ({\CH}) and ({\CL}) impliy that  
$$
I_{\beta}(k+1)-I_{\beta}(k)\sim\-\frac{\sqrt{\pi}}{4\sqrt{2}k^{3/2}},\ 
\ \ k\to\infty.\tag\CM
$$
Equations (\CK) and (\CM) determine the asymptotics of $D_{k+1}-D_k$,
and combining this with the asymptotics of $D_k$ given by (\CJ) we
obtain that the second term on the right hand side of (\CI) has
asymptotics $-1/(2k)$ as $k\to\infty$. The two terms on the right hand
side of (\CI) thus have a sum that is asymptotically
$-1/(2k)-1/(2k)=-1/k$, and Theorem~{\TAB} is proved. \quad \quad \qed

\enddocument